\RequirePackage{xcolor}
\documentclass[journal,twoside,web]{ieeecolor}
\usepackage{generic}
\usepackage{textcomp}
\usepackage{graphicx}      
\usepackage{siunitx}
\usepackage{pgfplots}
\pgfplotsset{compat=1.17}
\usepackage[acronym]{glossaries}
\usepackage{tikz}
\usetikzlibrary{decorations.pathmorphing}
\usetikzlibrary{shapes, backgrounds, patterns}
\usepackage{amsmath,verbatim,amssymb,amsfonts,amscd}
\usepackage{enumerate}
\usepackage{mathtools}
\usepackage{color}
\usepackage{mathrsfs}
\usepackage{cancel}
\usepackage{comment}
\usepackage{booktabs}

\newcommand{\dimM}{n}
\newcommand{\M}{\mathbb{R}^\dimM}
\newcommand{\SM}{\mathcal{S}}
\newcommand{\co}{z}

\newcommand{\dimu}{m}
\newcommand{\coords}{\co_1,\dots,\co_\dimM}

\newcommand{\tangentBundle}{T\M}

\newcommand{\partialDer}[2]{\frac{\partial#1}{\partial#2}}

\newcommand{\bmat}[1]{\begin{bmatrix}#1\end{bmatrix}}
\newcommand{\bfm}[1]{{\bf{#1}}}

\usepackage{chngcntr}
\usepackage{apptools}

\AtAppendix{\counterwithin{theorem}{section}}
\AtAppendix{\counterwithin{lemma}{section}}
\AtAppendix{\counterwithin{definition}{section}}
\AtAppendix{\counterwithin{proposition}{section}}
\AtAppendix{\counterwithin{remark}{section}}
\AtAppendix{\counterwithin{factnew}{section}}
\AtAppendix{\counterwithin{notation}{section}}

\newtheorem{theorem}{Theorem}
\newtheorem{lemma}{Lemma}
\newtheorem{prob}{Problem}

\newtheorem{definition}{Definition}
\newtheorem{proposition}{Proposition}
\newtheorem{remark}{Remark}

\newacronym{pde}{PDE}{partial differential equation}
\newacronym{nhim}{NHIM}{normally hyperbolic invariant manifold}
\newacronym{hjb}{HJB}{Hamilton-Jacobi-Bellman equation}
\newacronym{hje}{HJE}{Hamilton-Jacobi equation}
\newacronym{ode}{ODE}{Ordinary Differential Equation}
\newacronym{pre}{PRE}{periodic Riccati differential equation}

\usepackage{textcomp}
\def\BibTeX{{\rm B\kern-.05em{\sc i\kern-.025em b}\kern-.08em
    T\kern-.1667em\lower.7ex\hbox{E}\kern-.125emX}}
\begin{document}
\title{Optimal Stabilization of Periodic Orbits:\\ A Symplectic Geometry Approach} 

\author{Fabian Beck and Noboru Sakamoto
\thanks{The first author is with the School of Computation, Information and Technology, Technical University of Munich (TUM), 85748 Garching, Germany and also with the Institute~of~Robotics~and~Mechatronics, German~Aerospace~Center~(DLR), 82234 Wessling, Germany (e-mail:~fabian.beck@dlr.de).}
\thanks{The second author is with the Department of Mechanical Engineering and System Control, Nanzan University, Yamzato-cho 18, Showa-ku, Nagoya, 466-8673, Japan (e-mail: noboru.sakamoto@nanzan-u.ac.jp)}
\thanks{The work of the second author is partially funded by JSPS, Japan KAKENHI Grant Number JP22H01515 and by Nanzan University Pache Research Subsidy I-A-2 for 2024 academic year.}
}
\maketitle

\begin{abstract}                
In this contribution, the optimal stabilization problem of periodic orbits is studied via invariant manifold theory and symplectic geometry. The stable manifold theory for the optimal point stabilization case is generalized to the case of periodic orbit stabilization, where a \gls{nhim} plays the role of a hyperbolic equilibrium point. 
A sufficient condition for the existence of an \gls{nhim} of an extended Hamiltonian system is derived in terms of a periodic Riccati differential equation. It is shown that the problem of optimal orbit stabilization has a solution if a linearized periodic system is stabilizable and detectable. A moving orthogonal coordinate system is employed along the periodic orbit, which is a natural framework for orbital stabilization and linearization along the orbit. 
Two illustrative examples are presented: the first involves stabilizing a spring-mass oscillator at a target energy level, and the second addresses an orbit transfer problem for a satellite—a classic scenario in orbital mechanics. In both cases, we show that the proposed nonlinear feedback controller outperforms traditional linear control.
\end{abstract}


\section{Introduction}\label{sctn:introduciton}

For better or worse, periodic motions naturally arise in many branches of science and engineering. For instance, the repetitive firing of neurons—known to be related to neurophysiological disorders—can be described by limit cycles in systems of differential equations \cite{Hodgkin-1952,Hindmarsh198487}. In epidemiology, epidemic and endemic conditions are analyzed in terms of the stability or instability of certain limit cycles \cite{Kermack1932MathBio,Li1995Biosciences}. Heteroclinic connections between periodic orbits are relevant to capturing a comet into the gravitational field of a celestial body and to trajectory design for space missions \cite{Koon2000427}.
Furthermore, in our everyday lives, we rely on services provided by satellites, which must be stabilized to mission-dependent orbits \cite{Gurfil2016,Owis2009}.
Learning techniques for robotic locomotion in neural computation are developed using coupled oscillators in nonlinear dynamical systems \cite{Ijspeert2008642,Ijspeert2013328}. Finally, periodic gaits of biped robots are analyzed and designed via limit cycles \cite{Garcia1998281,Grizzle200151,Westervelt2003}.

From an engineering perspective—particularly in the control of mechatronic systems—the stabilization of periodic orbits has attracted considerable attention. Several approaches to the orbital stabilization problem exist.
One is based on virtual (holonomic) constraints \cite{Freidovich2008automatica,LaHera2013ieeeRobot,Kant2020sysconlett,Mohammadi2018automatica,Shiriaev2005tac,Shiriaev2006sysconlett,Shiriaev2008AnnRevCon,Shiriaev2014tac,Consolini2018}, while Herrera et al. \cite{Herrera2020ijc} propose a method that combines virtual constraints with robust stabilization techniques using a nonlinear $H^\infty$ method, as introduced in \cite{Orlov:2014:AHC}.
Another approach employs the immersion and invariance technique, proposed as a constructive design method for nonlinear and adaptive control in \cite{Astolfi2003ieeetac}. Relevant contributions along this line include \cite{Ortega2020rnc,Yi20208521,Romero2022ieeecst,Romero2024jrnc}.
In addition to the above approaches, Yi et al. \cite{Yi2020automatica} use the interconnection and damping assignment passivity-based control method originally proposed in \cite{Ortega2002automatica}, S{\ae}tre et al. \cite{Saetre2021rnc} apply sliding mode control, and S{\ae}tre and Shiriaev \cite{Saetre2023automatica} address the problem of inducing a stable heteroclinic orbit and performing point-to-point maneuvers in underactuated mechanical systems.

The present paper considers the optimal stabilization problem for periodic orbits.
The optimal point stabilization problem for nonlinear systems is one of the fundamental challenges in control theory and was first studied in \cite{Al'brekht:61:jamm,Lukes:69:sicon} using a Taylor series approach.
Van der Schaft demonstrated in \cite{vanderSchaft:91:syscon,vanderSchaft:92:ac} the advantages of a symplectic geometry-based approach for analyzing the \gls{hje} arising in the nonlinear $H^\infty$ control problem.

This paper extends the results of \cite{vanderSchaft:91:syscon,vanderSchaft:92:ac} from equilibrium points to more general steady-state behaviors—specifically, periodic orbits—by introducing several new tools.
The first is the use of a moving coordinate system (see \cite{Hale:73:ODE} for details and \cite{Banaszuk1995,Nielsen2008} for applications in geometric nonlinear control theory involving periodic orbits and submanifolds), which enables the measurement of deviations from the periodic orbit independently of its phase.
This approach excludes time-dependent tracking control problems for periodic solutions, such as those considered in \cite{HABAGUCHI2015215}.
The second key technique is the continuation method for periodic orbits, which allows for the systematic construction of an invariant manifold of the required dimension. This manifold serves as the foundation for applying the third and central tool in our framework: the theory of \glspl{nhim}.
These can be viewed as generalizations of hyperbolic equilibrium points and are characterized by invariant manifolds with attracting or repelling behavior.
We select an appropriate subset (a leaf) that plays the same role as the stable manifold in the case of optimal point stabilization in \cite{vanderSchaft:91:syscon,vanderSchaft:92:ac}.
Furthermore, the theory of \glspl{pre}, as studied in \cite{Kano1979:ijc,Kano1985:ieeetac}, plays a crucial role in linking the above concepts. In particular, it is instrumental in showing that the continued periodic orbits in the extended Hamiltonian system form an \gls{nhim}, where the relevant conditions are expressed through solutions to linear variational equations.

The main result (Theorem~\ref{thm:main_theorem}) is presented in a clear and concise form, and its conclusion follows naturally from the preceding analysis: a sufficient condition for the existence of an optimal orbital stabilizing controller is given in terms of stabilizability and detectability properties that must hold along the orbit.

The optimal stabilization problem for submanifolds, as studied in \cite{Montenbruck20172393}, is also applicable to the case of periodic orbits.
Here, we highlight the key differences between our approach and that of \cite{Montenbruck20172393}. 
In \cite{Montenbruck20172393}, it is shown that the submanifold stabilization problem has a local solution and that the associated feedback can be systematically computed using the algebraic Riccati equation, which corresponds to a standard infinite-horizon control problem.
This is achieved by utilizing a tubular neighborhood around the submanifold, which corresponds to the transverse coordinates in our setting.
Additionally, \cite{Montenbruck20172393} assumes that two matrices governing system controllability remain constant—an assumption noted as restrictive by the authors. 
In contrast, our approach employs the periodic Riccati equation and does not require this constraint.
Another key distinction is that \cite{Montenbruck20172393} does not account for tangential dynamics, which may hinder the asymptotic stability of the orbit. 

The present manuscript expands upon our conference contribution \cite{Beck2023IFAC}, incorporating more rigorous expositions of the \gls{nhim} construction using the continuation technique, as well as a second illustrative example.
The remainder of the paper is organized as follows. The optimal control problem addressed in this work is formulated in \S~\ref{sec:problem_statement}.
In addition, moving orthogonal coordinates are introduced to derive a system representation in which the state space is decomposed into tangential and transversal components.
In \S~\ref{sctn:continuation_NHIM}, we construct an invariant manifold for the Hamiltonian system associated with the \gls{hje} via continuation, and in \S~\ref{sctn:proof_NHIM}, we provide a detailed proof that this manifold is indeed an \gls{nhim}.
The local solution to the \gls{hje} is then constructed in \S~\ref{sctn:existence_OC} using laminations of the \gls{nhim}.
Two examples are presented in \S~\ref{sec:applications}: the first involves a mass-spring system stabilized at an orbit determined by a fixed energy level, while the second addresses an orbit transfer problem for a satellite—a canonical control challenge in orbital mechanics.
The paper concludes with a summary and final remarks in \S~\ref{sec:conclusion}, and the appendix contains technical details and prerequisites.
\section{Problem Statement}\label{sec:problem_statement}
\subsection{Base control system and its periodic orbit}
Let $\M$ be an $n$-dimensional euclidean state space with coordinates $\coords$ and consider a $C^r$ ($r\geqslant4$) dynamical control system 
\begin{equation}
    \dot{\co}=f(\co)+g(\co)u, \label{eqn:controlSys}
\end{equation}
where $g$ is a $\dimM\times\dimu$ matrix consisting of $\dimu$ $C^r$ vector fields as columns and $u\in\mathbb{R}^m$ represents the control inputs. 
Let $\SM\subset \M$ be a closed curve in $\mathbb{R}^n$ representing a periodic orbit of the unforced system (\ref{eqn:controlSys}) with $u=0$. 
By normalizing time, the period is set to 1. In Section~\ref{sec:tvc} a parametrization of $\SM$ is used;
\begin{equation}
    \SM = \{z\in\M | z=\gamma(x_0),\  0\leqslant x_0 \leqslant 1 \},\label{eqn:p-orbitParameterization}
\end{equation}
where $\gamma:\mathbb{R}/\mathbb{Z}\rightarrow\mathbb{R}^n$ is a 1-periodic $C^r$ curve, which is one-to-one on the interval $[0,1)$.

\subsection{Moving orthonormal system and cost functional}
\label{sec:tvc}
In this subsection, we introduce an orthonormal coordinate system in a neighborhood of the periodic orbit to enable a local and explicit analysis of the optimal stabilization problem.

We use the parametrization $\gamma(x_0)$ in (\ref{eqn:p-orbitParameterization}) of the original periodic orbit in the $z$-space, to construct a moving orthonormal frame along the periodic orbit, following the approach in \cite[Chapter~VI.1]{Hale:73:ODE}.
The moving orthonormal system about $\SM$ is constructed using
\begin{align*}
    e_0(x_0) = \left\lVert\partialDer{\gamma(x_0)}{x_0}\right\rVert^{-1}\partialDer{\gamma(x_0)}{x_0}
\end{align*}
together with $\dimM-1$ additional orthonormal vectors $e_1, \dots, e_{n-1}$ varying along $x_0$, i.e. $e_i:\mathbb{R}/\mathbb{Z} \rightarrow \mathbb{R}^n$ for $i=1,\dots, n-1$. 
With this, a coordinate transformation from $z$ to the new coordinate vector $\bfm{x} = [x_0, \bfm{x}_a] \in \mathbb{R} \times \mathbb{R}^{n-1}$ is defined by
\begin{align}
    z 
    = \gamma(x_0) + Z(x_0)\bfm{x}_a, \label{eqn:def_trafo_Hale}
\end{align}
where $x_0$ parametrizes the tangential direction with $0\leqslant x_0 \leqslant 1$, and $\bfm{x}_a \in \mathbb{R}^{n-1}$ represents the $(n-1)$ transverse components. The matrix $Z(x_0) = \begin{bmatrix} e_1(x_0) & \dots & e_{n-1}(x_0) \end{bmatrix}$ is a $\mathbb{R}^{n \times (n-1)}$-valued $C^{r-1}$ function that is 1-periodic in $x_0$.
Furthermore, it is known that there exists a constant $\rho > 0$, independent of $x_0$, such that the transformation \eqref{eqn:def_trafo_Hale} is valid for all $0\leqslant x_0\leqslant 1$ and $\|\bfm{x}_a\| < \rho$.
Applying the transformation \eqref{eqn:def_trafo_Hale} to \eqref{eqn:controlSys} yields
\begin{subequations}\label{eqn:controlSys_transformed}
\begin{align}
    \dot{x}_0 &= 1+f_0(x_0,\text{\bf x}_a)+g_0(x_0, \text{\bf x}_a)u,\\
    \dot{\bf x}_a &= A(x_0) \text{\bf x}_a + \text{\bf f}_a(\text{\bf x})+ g_a(x_0,\text{\bf x}_a)u,\label{eqn:controlSys_transformed_2nd}
\end{align}
\end{subequations}
where $f_0$, $g_0$, $A$, $\bfm{f}_a$ and $g_a$ are functions taking values in scalar, $m$-dimensional row vector, $(n-1)\times (n-1)$-matrix, $(n-1)$-dimensional vector and $(n-1)\times m$-matrix, respectively. 
The derivation of (\ref{eqn:controlSys_transformed}) from (\ref{eqn:controlSys}) can be shown as in \cite[page 217]{Hale:73:ODE}. 
In the same reference, Items (i)-(iii) of the following Lemma~\ref{lemma:periodic} are shown while Item (iv) follows from Taylor's formula and Item (i) (see also \cite{Whitney-1943}). 
\begin{lemma}\label{lemma:periodic}
\begin{enumerate}[(i)]
    \item 
$A$, $f_0$, $\bfm{f}_a$, $g_0$, $g_a$ are $C^{r-2}$ functions and period-1 in $x_0$ for all $\bfm{x}_a$. 
\item 
$|f_0(x_0,{\bf x}_a)|= \mathcal{O}(\|{\bf x}_a\|)$ as $\|{\bf x}_a\|\rightarrow 0$. 
\item 
${\bf f}_a(x_0,\bfm{0}_{n-1})=\bfm{0}_{n-1}$ and $\partialDer{\text{\bf f}_a}{\text{\bf x}_a}(x_0,\bfm{0}_{n-1})=\bfm{0}_{(n-1)\times (n-1)}$ for all $x_0\in\mathbb{R}$ and $\text{\bf x}_a\in\mathbb{R}^{\dimM-1}$.
\item
There is an $\mathbb{R}^{n-1}$-valued function $\Tilde{f}_0(x_0,\bfm{x}_a)$, which is $C^{r-3}$, such that $f_0(x_0,\bfm{x}_a)=\Tilde{f}_0(x_0,\bfm{x}_a)^\top \bfm{x}_a$.
\end{enumerate}
\end{lemma}

\vskip 1ex
In this paper, we limit ourselves to a penalty function on the state given by a quadratic function of the transverse coordinate $\bfm{x}_a$. 
Let $Q:\mathbb{R}\to \mathbb{R}^{(n-1)\times(n-1)}$ be a $C^r$ period-1 function of $x_0$ whose value is a positive semi-definite matrix, and consider a cost functional  
\begin{align}
    J(\bfm{x}(0),u)&=\int_0^\infty \frac{1}{2}\bfm{x}_a(t)^\top Q(x_0)\bfm{x}_a(t)+u^\top R u\,dt\label{eqn:cost_new_coordinates}
\end{align}
with the positive definite matrix $R\in\mathbb{R}^{m\times m}$ and the control set $L^2([0,\infty);\mathbb{R}^m)$, where $\bfm{x}(0)\in\mathbb{R}\times\mathbb{R}^{n-1}$ is the initial condition for (\ref{eqn:controlSys_transformed}). When the functional does not converge, we set $J=\infty$. 
The $\bfm{x}_a$-dynamics is said to be asymptotically stabilized by a control law if, for any neighborhood $U$ of the $x_0$-axis, there exists a neighborhood $U_0$ of the $x_0$-axis such that, for any initial condition in $U_0$, the corresponding closed-loop solution remains in $U$ for all $t\geqslant 0$ (Lyapunov stability), and if there further exists an open neighborhood $U_a$ of the $x_0$-axis such that, for any initial condition in $U_a$, the closed-loop solution $(x_0(t),\bfm{x}_a(t))$ satisfies $\dot x_0(t)\to 1$ and $\bfm{x}_a(t)\to 0$ as $t\to\infty$ (attractivity).
%
Now, the problem to be tackled in the paper is formulated as follows. 
\begin{prob}\label{prblm:main_problem}
Find, if it exists, a control law for (\ref{eqn:controlSys_transformed}) under which the 
$\bfm{x}_a$-dynamics (\ref{eqn:controlSys_transformed_2nd}) is asymptotically stabilized and, for each $\bfm{x}(0)$ in a neighborhood of $\bfm{x}_a=0$, the cost $J(\bfm{x}(0),u)$ is minimized over $u\in L^2([0,\infty);\mathbb{R}^m)$.
\end{prob}
\begin{remark}
Problem~\ref{prblm:main_problem} is a class of optimal control problems in the sense that closed-loop stability is required. This class of optimal control problems is called {\em stable regulator problem} and originated in \cite{Molinari:73:ac}. 
\end{remark}

\begin{remark}\label{rem:nominal_input}
The cost functional~\eqref{eqn:cost_new_coordinates} penalizes deviations from the periodic orbit independently of its phase.
For this, the periodic orbit must be a solution of the unforced system, that is, it must persist in the absence of the control input. Otherwise, the cost functional may not be convergent.
\end{remark}
\subsection{Solution approach}
To solve Problem~\ref{prblm:main_problem}, we take Bellman's sufficient condition approach and derive an associated \gls{hje}, a nonlinear \gls{pde} of the first order. 
Using the method of characteristics, an extended Hamiltonian system associated with the \gls{pde} is derived. We next show that the Hamiltonian system has a \gls{nhim}, where the theory of periodic Riccati equations is used. Using the theory of NHIMs, we show that the Hamiltonian system has a Lagrangian submanifold with attracting behaviors. Finally, using the integrability of the Lagrangian submanifold, the solution to the \gls{hje} is constructed. The derivative of the solution defines the control for Problem~\ref{prblm:main_problem}. 

\section{Continuation of trivial solution and a 2-dimensional invariant manifold}\label{sctn:continuation_NHIM}
In this section, we first derive an HJE (a \gls{pde}) associated with Problem~\ref{prblm:main_problem} and its characteristic equation. 

To apply dynamic programming, set
\begin{align*}
    H_d(\bfm{x},\bfm{p},u) & =  p_0 (1+f_0(\bfm{x})+g_0(\bfm{x})u)\\
    &\quad +\bfm{p}_a^\top \{A(x_0)\bfm{x}_a + \bfm{f}_a(\bfm{x})+ g_a(\bfm{x})u\}\\
    &\qquad +\frac{1}{2}\bfm{x}_a^\top Q(x_0)\bfm{x}_a + u^\top R u, 
\end{align*}
where $\mathbf{p} = [p_0, \mathbf{p}_a]^\top$, with $p_0 \in \mathbb{R}$ and $\mathbf{p}_a \in \mathbb{R}^{n-1}$ denoting the adjoint variables associated with $x_0$ and $\mathbf{x}_a$, respectively.
The optimality condition on $u$ (minimization of $H_d$ with respect to $u$) implies
\begin{align*}
    \partialDer{H_d}{u}=p_0 g_0(\bfm{x}) + \bfm{p}_a^\top g_a(\bfm{x})+ 2 u^\top R = 0.
\end{align*}
It follows that a minimizing $\bar{u}$ for $H_d$ is given by
\[
    \bar{u}(\bfm{x},\bfm{p}) = -\frac{1}{2}R^{-1}(g_0(\bfm{x})^\top p_0 + g_a(\bfm{x})^\top \bfm{p}_a).
\]
We introduce the shorthand notation
\[
    G(\bfm{x},\bfm{p}) = g_0(\bfm{x})^\top p_0 + g_a(\bfm{x})^\top \bfm{p}_a,
\]
to simplify the expression in the following. 
Using this, we define 
\begin{align*}
    H(\bfm{x},\bfm{p})  & = H_d(\bfm{x},\bfm{p},\bar{u}(\bfm{x},\bfm{p}))\\
    &=
    p_0+ p_0 f_0(\bfm{x})
    -\frac{1}{4}G(\bfm{x},\bfm{p})^\top R^{-1} G(\bfm{x},\bfm{p})\\
     &\quad +\bfm{p}_a^\top (A(x_0)\bfm{x}_a + \bfm{f}_a(\bfm{x}))+\frac{1}{2}\bfm{x}_a^\top Q(x_0)\bfm{x}_a.
\end{align*}

By setting $H(\bfm{x},\bfm{p})$ to zero we obtain the \gls{hje} in form of a \gls{pde}
\begin{align}
    H\left(\bfm{x},\frac{\partial V}{\partial \bfm{x}}\right) = 0,\label{eqn:HJE_newCordinate}
\end{align}
where we wish to find the differentiable function $V(\bfm{x})$ with certain properties to solve Problem~\ref{prblm:main_problem}. 
The corresponding characteristic equations (see, e.g., \cite[Appendix~7]{Libermann:87:SGAM}) for this \gls{pde} form a Hamiltonian system, obtained by extending the state with $n$ adjoint variables. This extended system is given by
\begin{subequations}\label{eqn:hamiltoniansys}
\begin{align}
    \dot{x}_0 &= \!1\!+\!f_0\!({\bf x})\!-\!\frac{1}{2}g_0\!({\bf x} )R^{-1}G({\bf x},{\bf p}),\\
    \dot{{\bf x}}_a &= \!A(x_0){\bf x}_a\!+\!{\bf f}_a({\bf x})\!-\frac{1}{2}{g}_a({\bf x})R^{-1}G({\bf x},{\bf p}),\\
    \dot{p}_0 &= -p_0\partialDer{f_0}{x_0}(\bfm{x})+\frac{1}{4}\partialDer{}{x_0}(G(\bfm{x},\bfm{p})^\top R^{-1}G(\bfm{x},\bfm{p})) \nonumber\\
    &- \bfm{p}_a^\top \frac{dA}{dx_0}(x_0)\bfm{x}_a - \bfm{p}_a^\top \partialDer{\bfm{f}_a}{x_0}(\bfm{x})- \frac{1}{2}\bfm{x}_a^\top \frac{dQ}{dx_0}\bfm{x}_a, \label{eqn:hamiltoniansys_p0}\\
    \dot{\bfm{p}}_a &= -p_0 \partialDer{f_0}{\bfm{x}_a}(\bfm{x})^\top\!\!+\!\frac{1}{4}\!\left(\!\partialDer{}{\bfm{x}_a}(G(\bfm{x},\bfm{p})^\top\!R^{-1}G(\bfm{x},\bfm{p}))\!\right)^\top \nonumber\\
    &-A(x_0)^\top \bfm{p}_a - \partialDer{\bfm{f}_a}{\bfm{x}_a}(\bfm{x})^\top \bfm{p}_a - Q(x_0)\bfm{x}_a.\label{eqn:hamiltoniansys_pa}
\end{align}
\end{subequations}
The right-side of \eqref{eqn:hamiltoniansys} is a Hamiltonian vector field with Hamiltonian $H(\bfm{x},\bfm{p})$ and this vector field is denoted as $X_H(\bfm{x},\bfm{p})$. 
Also, in what follows, $\Phi_H(t,(\bfm{x},\bfm{p}))$ denotes its flow starting from $(\bfm{x},\bfm{p})$ at $t=0$. 
The following Lemma demonstrates that there is an invariant manifold for (\ref{eqn:hamiltoniansys}) in $\mathbb{R}^{2n}$ corresponding to the periodic orbit $\gamma$ in the $\bfm{x}$-space  or to the trivial solution $(t,\bfm{0}_{2n-1})$ of (\ref{eqn:hamiltoniansys}). 
\begin{lemma}\label{lemma:invariance}
The manifold 
\[
\Gamma_0=\{(\bfm{x},\bfm{p})\,|\,x_0\in \mathbb{R},\bfm{x}_a=\bfm{0}_{n-1},p_0=0,\bfm{p}_a=\bfm{0}_{n-1}\}
\]
is invariant under the flow of $X_H$.
\end{lemma}

\begin{proof}
Since $\mathbf{p} = 0$ and $\bfm{x}_a=0$ on $\Gamma_0$, it follows from 
Lemma~\ref{lemma:periodic} property (iii) that \[
f_0(x_0,\bfm{0}_{n-1})=0,G(x_0,\bfm{0}_{2n-1})=\bfm{0}_m,
\bfm{f}_a(x_0,\bfm{0}_{n-1})=\bfm{0}_{n-1}
\]
for all $x_0$.
Thus, $\dot{\mathbf{p}} = 0$ and $\dot{x}_0 = 1$ as well as $\dot{\bfm{x}}_a = \bfm{0}_{n-1}$, which establishes the invariance $\Gamma_0$.
\end{proof}

The dynamics of \eqref{eqn:hamiltoniansys} is linearized along the trivial solution $(t,\bfm{0}_{2n-1})$ and is characterized by
\begin{align}
&	\mathscr{H}(t)= \begin{bmatrix}
		0 & \tilde{f}_0(t,\bfm{0}_{n-1})^\top &W_{00}(t) & W_{0a}(t)\\
		\bfm{0}_{n-1} & {A}(t) & W_{a0}(t)  &-\bar{R}(t)\\					
		0 & \bfm{0}_{n-1}^\top & 0 & \bfm{0}_{n-1}^\top\\
		\bfm{0}_{n-1} & - Q(t)& -\tilde f_0(t,\bfm{0}_{n-1}) & -{A}(t)^\top
		\end{bmatrix},\label{eqn:lin_hamiltoniansys}\\
\intertext{where}
    &W_{ij}(t)=-\frac{1}{2} g_i(t,\bfm{0}_{n-1}) R^{-1} g_j(t,\bfm{0}_{n-1})^\top, \ i,j=0,a \nonumber\\
    &\bar{R}(t)=-W_{aa}(t)=\frac{1}{2}g_a(t,\bfm{0}_{n-1}) R^{-1}g_a(t,\bfm{0}_{n-1})^\top\nonumber
\end{align}
and $\Tilde{f}_0(x_0,\bfm{x}_a)$ is defined as in Lemma~\ref{lemma:periodic}. 
Namely, the variational equation of (\ref{eqn:hamiltoniansys}) along the trivial solution is
\[
\frac{d}{dt}\bmat{\bar{x}_0\\ \bar{\bfm{x}}_a\\ \bar{p}_0\\ \bar{\bfm{p}}_a}
=\mathscr{H}(t)
\bmat{\bar{x}_0\\ \bar{\bfm{x}}_a\\ \bar{p}_0\\ \bar{\bfm{p}}_a}, 
\]
where $\bar{x}_0$, $\bar{\bfm{x}}_a$, $\bar{p}_0$, and $\bar{\bfm{p}}_a$ denote the variational quantities representing deviations from the trivial solution, and correspond to $x_0$, $\bfm{x}_a$, $p_0$, and $\bfm{p}_a$, respectively.

The objective here is to show that the invariant manifold $\Gamma_0$ containing the trivial solution of (\ref{eqn:hamiltoniansys}) can be extended by continuation for nonzero Hamiltonian values from which an \gls{nhim} is constructed. 
To this end, we will derive conditions such that the monodromy matrix $D\Phi_H(1,\bfm{0}_{2n})$ for the period-1 variational equation with (\ref{eqn:lin_hamiltoniansys}) has eigenvalue $1$ with multiplicity $2$. 
It is convenient to change the order of variables from $(\bar{x}_0, \bar{\bfm{x}}_a, \bar{p}_0, \bar{\bfm{p}}_a)$ to $(\bar{x}_0, \bar{p}_0, \bar{\bfm{x}}_a, \bar{\bfm{p}}_a)$ (hereafter, we call this coordinate change `reordering'). 
Correspondingly, the transformed vector field will be denoted as $\Tilde{X}_H$ and its flow will be $\Tilde{\Phi}_H$. Then (\ref{eqn:lin_hamiltoniansys}) becomes
\begin{equation}
\Tilde{\mathscr{H}}(t)= \begin{bmatrix}
    0 &W_{00}(t)  &\tilde{f}_0(t,\bfm{0}_{n-1})^\top & W_{0a}(t)\\
  0 & 0  & \bfm{0}_{n-1}^\top & \bfm{0}_{n-1}^\top\\
    \bfm{0}_{n-1} & W_{a0}(t) & {A}(t)  &-\bar{R}(t)\\					
    \bfm{0}_{n-1} & -\tilde f_0(t,\bfm{0}_{n-1})& - Q(t) & -{A}(t)^\top
		\end{bmatrix}.\label{eqn:lin_hamiltoniansys_new}
\end{equation}
We notice that 
\[
\mathrm{Ham}(t)=\bmat{
 {A}(t)  &-\bar{R}(t)\\					
 - Q(t) & -{A}(t)^\top
    },
\]
which plays an important role hereafter, appears in the right-bottom block. 
Now, the following observation will be useful. 
\begin{lemma}\label{lemma:form_of_sol_for_H}
The fundamental matrix ${M}_{\Tilde{\mathscr{H}}}(t,0)$ of (\ref{eqn:lin_hamiltoniansys_new}) has the form of
\begin{equation}
\bmat{1& c(t)& \bfm{k}(t)^\top\\
0& 1& \bfm{0}_{2n-2}^\top\\
\bfm{0}_{2n-2}& \bfm{h}(t)& M_{\mathrm{Ham}}(t,0)
},\label{eqn:FundamentalMatrix}
\end{equation}
where $c(t)$ is a continuous scalar function of $t$,  $\bfm{k}(t)$ and $\bfm{h}(t)$ are $\mathbb{R}^{2n-2}$-valued continuous functions of $t$ and $M_{\mathrm{Ham}}(t,0)$ is the fundamental matrix for $\mathrm{Ham}(t)$. 
\end{lemma}

\vskip 1ex
Lemma~\ref{lemma:form_of_sol_for_H} is shown by using twice the fact that the fundamental matrix of a linear differential equation with a block-triangular periodic matrix has a block-triangular fundamental matrix. Thus, it suffices to show that the monodromy matrix $M_{\mathrm{Ham}}(1,0)$ of
$\mathrm{Ham}(t)$ has no eigenvalue $1$. $\mathrm{Ham}(t)$ has a well-known relation with a periodic differential Riccati equation
\begin{multline}
    \dot{P}(t) + A(t)^\top P(t) + P(t)A(t) \\
     -P(t)\bar{R}(t)P(t) +Q(t)= 0.\label{eqn:p1Riccati}
\end{multline}
A period-1 solution $P(t)$ is said to be a {\em stabilizing solution} for (\ref{eqn:p1Riccati}) if $A(t)-R(t)P(t)$ is asymptotically stable. See Appendix~\ref{sec_app:linearPeriodicSystems} for the theory of linear periodic systems and periodic Riccati differential equations. 
\begin{proposition}\label{prop:pRic_eigenvalues}
Suppose that the periodic Riccati equation (\ref{eqn:p1Riccati}) has a period-1 stabilizing solution $P_s(t)$. Then, the monodromy matrix of $\mathrm{Ham}(t)$ has $n-1$ eigenvalues inside the unit circle and $n-1$ eigenvalues outside of the unit circle. Moreover, $D\Phi_H(1,\bfm{0}_{2n})$ has eigenvalue 1 with algebraic multiplicity 2. 
\end{proposition}
\begin{proof}
We first consider a period-1 linear differential equation associated with $\mathrm{Ham}(t)$
\begin{align}
    \frac{d}{dt}\bmat{\bar{\bfm{x}}_a\\ \bar{\bfm{p}}_a} &= \bmat{{A}(t)&-\bar{R}(t) \\- Q(t)& -{A}(t)^\top}\bmat{\bar{\bfm{x}}_a\\ \bar{\bfm{p}}_a}.\label{eqn:linearTransverse}
\end{align}
Note that it holds that 
\begin{equation}
    J\mathrm{Ham}(t)+\mathrm{Ham}(t)^\top J=0,\label{eqn:simplectic_diff}
\end{equation}
where $J=\left[\begin{smallmatrix}0&I\\-I&0\end{smallmatrix}\right]$ and $I$ is the identity matrix of $n-1$ dimension. 
Let $\bar{x}_j(t)$, $j=1,\ldots,n-1$, be independent solutions to 
\begin{equation}
\frac{d}{dt}\bar{\bfm{x}}_a=({A}(t)-\bar{R}(t)P_s(t))\bar{\bfm{x}}_a.\label{eqn:closed_loop}
\end{equation}
Since $P_s(t)$ is assumed to be a stabilizing solution, it follows that $\bar{x}_j(t) \to 0$ as $t \to \infty$ for $j = 1, \ldots, n-1$.
Also, any solution $\bar{\bfm{x}}_a(t)$ of (\ref{eqn:closed_loop}) satisfies
\begin{align*}
&\frac{d}{dt}\bmat{\bar{\bfm{x}}_a(t)\\P_s(t)\bar{\bfm{x}}_a(t)} \\&=\bmat{({A}(t)-\bar{R}(t)P_s(t))\bar{\bfm{x}}_a(t)\\
\dot{P}_s \bar{\bfm{x}}_a+P_s\dot{\bar{\bfm{x}}}_a}\\
&=\bmat{({A}-\bar{R}P_s)\bar{\bfm{x}}_a\\
-(P_s{A}+{A}^\top P_s-P_s\bar{R}P+Q)\bar{\bfm{x}}_a+P_s({A}-\bar{R}P_s)\bar{\bfm{x}}_a}\\
&=\bmat{{A}(t)-\bar{R}(t)P_s(t)\\-Q(t)-{A}(t)^\top P_s(t)}\bar{\bfm{x}}_a\\
&=\bmat{{A}(t)&-\bar{R}(t)\\ -Q(t)& -{A}(t)^\top}\bmat{\bar{\bfm{x}}_a(t)\\P_s(t)\bar{\bfm{x}}_a(t)},
\end{align*}
where we have used (\ref{eqn:p1Riccati}). This shows that $\left[\begin{smallmatrix}
\bar{x}_j(t)\\P_s(t)\bar{x}_j(t)\end{smallmatrix}\right]$, $j=1,\ldots,n-1$, are $n-1$ independent solutions of (\ref{eqn:linearTransverse}) that converge to $0$ as $t\to\infty$ since $P_s(t)$ is a period-1 matrix. 

Let $M_{\mathrm{Ham}}(t,s)$ be the transition matrix for the linear differential equation (\ref{eqn:linearTransverse}) and set $\zeta_j=\left[\begin{smallmatrix}\bar{x}_j(0)\\P_s(0)\bar{x}_j(0)\end{smallmatrix}\right]$, $j=1,\ldots,n-1$. Then we have 
\[
\bmat{\bar{x}_j(t)\\P_s(t)\bar{x}_j(t)}=M_{\mathrm{Ham}}(t,0)\zeta_j, \ j=1,\ldots,n-1.
\]
The fact we have just shown means that for $k\in\mathbb{N}$
\[
M_{\mathrm{Ham}}(k,0)\zeta_j\to0\ \text{as }k\to\infty, \ j=1,\ldots,n-1.
\]
From the periodicity, we have $M_{\mathrm{Ham}}(k,0)=M_{\mathrm{Ham}}(1,0)^k$. Now, we can show from (\ref{eqn:simplectic_diff}) that if $\lambda\in\mathbb{C}$ is an eigenvalue of $M_{\mathrm{Ham}}(1,0)$, so is $1/\lambda$ (detail is omitted). Thus we conclude that $\zeta_1,\ldots,\zeta_{n-1}$ belong to the generalized eigenspace corresponding to the eigenvalues of $M_{\mathrm{Ham}}(1,0)$ located inside the unit circle. Because eigenvalues of $M_{\mathrm{Ham}}(1,0)$ are located symmetrically with respect to the unit circle, $M_{\mathrm{Ham}}(1,0)$ has $n-1$ eigenvalues outside the unit circle. 
The last statement is proved from Lemma~\ref{lemma:form_of_sol_for_H}. 
\end{proof}
\begin{proposition}\label{prop:continuation_for_Hamsys}
Suppose that the periodic Riccati equation (\ref{eqn:p1Riccati}) has a period-1 stabilizing solution $P_s(t)$. Then, there exist an $\epsilon_1>0$ and $C^{r-2}$ univariate functions $\tau:[-\epsilon_1,\epsilon_1]\to\mathbb{R}$, $\bfm{a}:[-\epsilon_1,\epsilon_1]\to\mathbb{R}^{2n-1}$ such that $\tau(0)=1$, $\bfm{a}(0)=\bfm{0}_{2n-1}$ and the solution $(x_0(t),\bfm{x}_a(t),\bfm{p}(t))$ of (\ref{eqn:hamiltoniansys}) starting from $(x_0(0),\bfm{a}(e))$ at $t=0$ satisfy for all $e\in [-\epsilon_1,\epsilon_1]$ the following conditions: 
\begin{enumerate}[(i)]
    \item $x_0(m\tau(e))=x_0(0)+m$ for $m\in\mathbb{N}$,
    \item $\bfm{x}_a(t)$  and $\bfm{p}(t)$ are period-$\tau(e)$ functions.
\end{enumerate}
Furthermore, 
\begin{multline*}
\Gamma(e)=\{(x_0(t),\bfm{x}_a(t),\bfm{p}(t))\,|\, x_0(0)=0, \\ 
    (\bfm{x}_a(0),\bfm{p}(0))=\bfm{a}(e),t\in\mathbb{R}  \}
\end{multline*}
is a family of $C^{r-2}$ invariant manifolds for (\ref{eqn:hamiltoniansys}) such that 
$\Gamma(0)=\Gamma_0$, $\Gamma(e)\subset H^{-1}(e)$. 
\end{proposition}
\begin{proof}
For the sake of notation simplicity, let us introduce 
\[y_0=x_0,\ \bfm{y}_a=(\bfm{x}_a,\bfm{p})\in\mathbb{R}^{2n-1},\ \bfm{y}=(y_0,\bfm{y}_a)
\]
and new functions $h_0(y_0.\bfm{y}_a)$, $\bfm{h}_a(y_0.\bfm{y}_a)$, to rewrite (\ref{eqn:hamiltoniansys}) as
\begin{gather*}
\dot{y}_0=1+h_0(y_0,\bfm{y}_a),\quad \dot{\bfm{y}}_a=\bfm{h}_a(y_0.\bfm{y}_a),\\
\intertext{where, }
\begin{aligned}
&\text{the right-side of (\ref{eqn:hamiltoniansys}a)}=1+h_0(y_0.\bfm{y}_a),\\
&\text{the right-sides of (\ref{eqn:hamiltoniansys}b)$\sim$(\ref{eqn:hamiltoniansys}d)}=\bfm{h}_a(y_0.\bfm{y}_a).
\end{aligned}
\end{gather*}
We may also rewrite (\ref{eqn:hamiltoniansys}) as $\dot{\bfm{y}}=\bfm{h}(\bfm{y})$, where $\bfm{h}(\bfm{y})=[1+h_0(\bfm{y})\ \bfm{h}_a(\bfm{y})^\top] ^\top$. 
The following properties can be shown from Lemma~\ref{lemma:periodic}. 
\begin{enumerate}[(a)]
    \item $h_0(y_0,\bfm{y}_a)$, $\bfm{h}_a(y_0.\bfm{y}_a)$ are of class $C^{r-2}$ and period-1 functions in $y_0$ for all $\bfm{y}_a$. 
    \item $h_0(y_0,\bfm{0}_{2n-1})=0$, $\bfm{h}_a(y_0,\bfm{0}_{2n-1})=\bfm{0}_{2n-1}$ for $y_0\in\mathbb{R}$.
    \item $H(y_0,\bfm{0}_{2n-1})=0$ for $y_0\in\mathbb{R}$.
    \item $dH\neq 0$ on $\Gamma_0=\{(y_0,\bfm{0}_{2n-1})\,|\,y_0\in\mathbb{R}\}$. 
\end{enumerate}
\newcounter{temp}
\setcounter{temp}{\value{enumi}}
From the conservation of $H$ along (\ref{eqn:hamiltoniansys}) and the above periodicity properties, we have
\begin{enumerate}[(a)]
\setcounter{enumi}{\value{temp}}
    \item $dH(\bfm{y})\bfm{h}(\bfm{y})=0$ for $\bfm{y}\in\mathbb{R}^{2n}$. 
    \item $D\Phi_H(1,\bfm{y}^\mathrm{o})\bfm{h}(\bfm{y}^\mathrm{o})=\bfm{h}(\bfm{y}^\mathrm{o})$ for $\bfm{y}^\mathrm{o}\in\Gamma_0$. 
    \item $dH(\bfm{y}^\mathrm{o})D\Phi_H(1,\bfm{y}^\mathrm{o})=dH(\bfm{y}^\mathrm{o})$ for $\bfm{y}^\mathrm{o}\in\Gamma_0$. 
\end{enumerate}

Now we employ a $C^{r-2}$ coordinate change from $\bfm{y}_a$ to $(H,\bfm{y}_b)$, $(y_0,\Psi(\bfm{y}_a))=(y_0,H,\bfm{y}_b)$, which is defined in a neighborhood of $\Gamma_0$, with $\bfm{y}_b\in\mathbb{R}^{2n-2}$ in such a way that   
\[
dx_0\cdot\bfm{y}_b=0,\ dH\cdot\bfm{y}_b=0.
\]
In the new coordinates, a solution of (\ref{eqn:hamiltoniansys}) starting from $(c,e,\bfm{y}_b)$ at $t=0$ is written as
\[
\left(y_0(t,(c,e,\bfm{y}_b)),\ e,\ K(t,(c,e,\bfm{y}_b))\right), 
\]
where $K$ is a $C^{r-2}$ function from $\mathbb{R}\times\mathbb{R}\times\mathbb{R}\times\mathbb{R}^{2n-2}$ to $\mathbb{R}^{2n-2}$. Noting that $\frac{\partial}{\partial t} y_0(0,\bfm{0}_{2n})=1\neq0$, from the implicit function theorem, there exists a $C^{r-2}$ function $\tilde{\tau}(e,\bfm{y}_b)$ defined for sufficiently small $|e|$ and $\|\bfm{y}_b\|$ such that
\[
\Tilde{\tau}(0,\bfm{0}_{N-2})=1,\ 
y_0(\Tilde{\tau}(e,\bfm{y}_b),(0,e,\bfm{y}_b))=1
\]
for all $e$, $\bfm{y}_b$ for which $\Tilde{\tau}$ is defined. 
One can show, using for instance Lemma~\ref{lemma:linalg_app} below, that the monodromy matrix $D\Phi_H(1,\bfm{0}_{2n})$ is written in the new coordinates as
 \[
 \left[
\begin{array}{cc|ccc}
1& *& * & \cdots& *\\
0& 1& 0& \cdots&0\\\hline
0& *& & & \\
\vdots & \vdots & & \frac{\partial K}{\partial \bfm{y}_b}(1,\bfm{0}_{2n}) & \\
0& * & & &
\end{array}\right], 
\]
and that $D\Phi_H(1,\bfm{0}_{2n})$ has eigenvalue 1 with algebraic multiplicity larger than or equal to 2. 
Now, we apply Proposition~\ref{prop:pRic_eigenvalues} and conclude that 
$\frac{\partial K}{\partial \bfm{y}_b}(0,\bfm{0}_{2n})$ does not have eigenvalues at 1. This shows that, by the implicit function theorem, the periodicity requirement for the 3rd$\sim\!2n\,$th-components; 
\[
K(\Tilde{\tau}(e,\bfm{y}_b),(0,e,\bfm{y}_b))=\bfm{y}_b
\]
has a $C^{r-2}$ solution $\bfm{y}_b=\eta_b(e)$ with $\eta_b(0)=\bfm{0}_{2n-2}$ for $|e|\leqslant\epsilon_1$, where $\epsilon_1>0$ is taken in such a way that the 1st and 2nd applications of the implicit function theorem are guaranteed. By setting 
\begin{gather*}
\tau(e)=\Tilde{\tau}(e,\eta_b(e)),\\
\bfm{a}(e)=\Psi^{-1}(e,\eta_b(e)),
\end{gather*}
we have the following evolution: 
\[
\bmat{0\\ \bfm{a}(e)}\xrightarrow{\tau(e)}\bmat{1\\ \bfm{a}(e)}\xrightarrow{\tau(e)}\bmat{2\\ \bfm{a}(e)}\xrightarrow{\tau(e)}\cdots
\]
which shows Items (i) and (ii) as well as the last statement. 
\end{proof}
\begin{lemma}\label{lemma:linalg_app}
Let $A$ be an $n\times n$ matrix and assume that there exist two nonzero vectors $u$, $v$ satisfying 
\[
Au=u,\quad v^\top A=v^\top,\quad v^\top u=0,
\]
then $A$ has eigenvalue $1$ with algebraic multiplicity larger than or equal to $2$. 
\end{lemma}
\begin{proof}
The Lemma is proved using the Gram-Schmidt procedure. 
\end{proof}
\begin{remark}
The solution corresponding to $(0,\bfm{a}(e))$ is called {\em a continuation} of the trivial solution. 
The proof of Proposition~\ref{prop:continuation_for_Hamsys} uses the technique in \cite[p.~160]{Meyer2017} or \cite[p.~436]{Chicone2006} which is called {\em the Lyapunov-Schmidt reduction} (see, also \cite[p.~496]{Abraham:79:FOM}). In the present case, however, the trivial orbit and its continuation are not periodic orbits and some modifications need to be made. 
\end{remark}
\vskip 1ex
Since all the functions in \eqref{eqn:hamiltoniansys} have period $1$ in $x_0$, using above (i)-(ii) in Proposition~\ref{prop:continuation_for_Hamsys}, we take a quotient space 
$\mathbb{R}/\mathbb{Z}$ for $x_0$-coordinate to get a family of rings $\Gamma(e)/\mathbb{Z}$. Now we have a 2-dimensional compact invariant $C^{r-1}$ manifold 
\[\mathscr{M}(\epsilon_1)=\bigcup_{|e|\leqslant\epsilon_1}\Gamma(e)/\mathbb{Z}\]
for \eqref{eqn:hamiltoniansys} (see Fig.~\ref{fig:invariantCompactManifold}). 
\begin{remark}
In the next section, we prove that $\mathscr{M}(\epsilon)$ (for small $\epsilon$) is an \gls{nhim} (see Appendix~\ref{sctn:nhim_appx} for definition). It has a stable manifold and a subset (leaf) of the stable manifold will be the graph of a derivative of the desired solution to HJE (\ref{eqn:HJE_newCordinate}). We, however, note that $\Gamma_0$ itself cannot be an \gls{nhim} because dimension counts do not agree (see (\ref{eqn:split_appdx}) in the Appendix~\ref{sctn:nhim_appx} and note that the dimensions of $N_x^s$ and $N_x^u$ are $n-1$ in our case). This is why we employed the continuation. We also note that the compactness required in the \gls{nhim} theory is satisfied from the periodicity (see Item (i)-(ii) in Proposition~\ref{prop:continuation_for_Hamsys}). We can either take a quotient above or work directly in the $(x_0,\bfm{x}_a)$ coordinates carefully requiring periodicity in $x_0$ for functions that additionally appear. The latter approach will be taken in the solution construction for HJE (\ref{eqn:HJE_newCordinate}) (see the proof of Proposition~\ref{prop:lagrange_SM_lammination}).
\end{remark}
\begin{figure}
    \centering
\begin{tikzpicture}
\begin{axis}[
    hide axis,
    view={40}{40}
]
\addplot3 [
    surf, shader=faceted interp,
    point meta=x,
    colormap/greenyellow,
    samples=40,
    samples y=11,
    z buffer=sort,
    domain=0:360,
    y domain=-0.5:0.5
] (
    {(1+0.5*y*cos(x)))*cos(x)},
    {(1+0.5*y*cos(x)))*sin(x)},
    {0.5*y*sin(x)});

\addplot3 [
    samples=50,
    domain=-120:195,
    samples y=0,
    thick
] (
    {cos(x)},
    {sin(x)},
    {0});
\draw[->] (xyz cs:y=-1.5) -- (xyz cs:y=-1.2,z=0.065) ;
\node at (xyz cs:x=0.1,y=-2,z=0.1)     {$e = 0$};
\draw[->] (xyz cs:y=-1.02,z=-0.2) -- (xyz cs:y=-1.05,z=0.2) ;
\node at (xyz cs:y=-0.8,z=0.1)     {$e \uparrow $};
\end{axis}

\end{tikzpicture}
    \caption{The invariant manifold $\mathscr{M}(\epsilon)$ and the orbit given by $\Gamma_0$ (black line).}
    \label{fig:invariantCompactManifold}
\end{figure}

\section{Proof that $\mathscr{M}$ is an NHIM}\label{sctn:proof_NHIM}
Let $\Tilde{\mathscr{H}}(t)$ in (\ref{eqn:lin_hamiltoniansys_new}) be partitioned as 
\[
\Tilde{\mathscr{H}}(t)=\bmat{H_{11}(t) & H_{12}(t)\\
    H_{21}(t)& \mathrm{Ham}(t)},
\]
where $H_{11}$, $H_{12}$ and $H_{21}$ are $2\times2$, $2\times (2n-2)$ and $(2n-2)\times2$ matrices, respectively and introduce corresponding $2$ and $2n-2$ variational states 
\[
\xi=\bmat{\bar{x}_0\\ \bar{p}_0},\ \eta=\bmat{\bar{\bfm{x}}_a\\\bar{\bfm{p}}_a}
\]
so that the variational equation of (\ref{eqn:hamiltoniansys}) along the trivial solution ($e=0$) is equivalently (using the new order of variables) written as
\begin{align}
\frac{d}{dt}\bmat{\xi\\\eta}
=\bmat{H_{11}(t) & H_{12}(t)\\
    H_{21}(t)& \mathrm{Ham}(t)}\bmat{\xi\\\eta}.  \label{eqn:linearized_dynamics}
\end{align}
For nonzero $e\in [-\epsilon_1,\epsilon_1]$, let us write the variational equation along the solution evolving in $\Gamma(e)$ as 
\begin{align}
\frac{d}{dt}\bmat{\xi^e\\\eta^e}
=\bmat{H_{11}^e(t) & H_{12}^e(t)\\
    H_{21}^e(t)& H_{22}^e(t)}\bmat{\xi^e\\\eta^e}. \label{eqn:linearized_dynamics_nonzero_e}
\end{align}
Note first that $H_{ij}^e(t)$ ($i,j=1,2$) are period-$\tau(e)$ matrices and  second that, from $C^{r-1}$ smoothness, it follows that 
\begin{gather}
H_{11}^e(t)\to H_{11}(t), \ H_{12}^e(t)\to H_{12}(t), \ H_{21}^e(t)\to H_{21}(t),\nonumber\\
H_{22}^e(t)\to \mathrm{Ham}(t) \label{eqn:convergence_He2H}
\end{gather}
for $t\in\mathbb{R}$ as $e\to 0$. 

Let $\Phi_H(t,(x_0,\bfm{a}(e)))$ denote the solution of \eqref{eqn:hamiltoniansys} starting from $(x_0,\bfm{a}(e))\in \mathscr{M}(\epsilon_1)$ at $t=0$, where $\bfm{a}(e)$ is as defined in Proposition~\ref{prop:continuation_for_Hamsys}, and then $D\Phi_H(\tau(e),(x_0,\bfm{a}(e)))$ is the monodromy matrix of the corresponding variational equation along $\Phi_H(t,(x_0,\bfm{a}(e)))$. Let also $\Tilde{\Phi}_H(t,(x_0,\bfm{a}(e)))$ denote the solution of reordered \eqref{eqn:hamiltoniansys} starting from $(x_0,\bfm{a}(e))\in \mathscr{M}(\epsilon_1)$ at $t=0$. Then, $D\Tilde{\Phi}_H(\tau(e),(x_0,\bfm{a}(e)))$ is the monodromy matrix for (\ref{eqn:linearized_dynamics_nonzero_e}). $D{\Phi}_H(\tau(e),(x_0,\bfm{a}(e)))$ and $D\Tilde{\Phi}_H(\tau(e),(x_0,\bfm{a}(e)))$ are similar to each other via the Jacobian matrix of the reordering and we show the properties of the former using the latter. 

Now, we have the following theorem. 
\begin{theorem}\label{thm:mainestimates}
Assume that the periodic Riccati equation \eqref{eqn:p1Riccati} associated with $\mathrm{Ham}(t)$ has a period-1 stabilizing solution. Then, there exists a positive $\epsilon_2<\epsilon_1$ such that for $\mathscr{M}(\epsilon_2)$, the following hold. 
\begin{enumerate}[(i)]
    \item For $(x_0,\bfm{a}(e))\in\mathscr{M}(\epsilon_2)$, $0\leqslant x_0<1$, $0\leqslant e<\epsilon_2$, $\mathbb{R}^{2n}$ we have an $D{\Phi}_H(\tau(e),(x_0,\bfm{a}(e)))$-invariant splitting
    \[
\mathbb{R}^{2n}=T_{(x_0,\bfm{a}(e))}\mathscr{M}(\epsilon_2)\oplus N^{s}_{(x_0,\bfm{a}(e))} \oplus N^{u}_{(x_0,\bfm{a}(e))},
    \]
where $T_{(x_0,\bfm{a}(e))}\mathscr{M}(\epsilon_2)$, $N^{s}_{(x_0,\bfm{a}(e))}$ and $N^{u}_{(x_0,\bfm{a}(e))}$ are $2$, $n-1$ and $n-1$ dimensional subspaces, respectively, which are all invariant under $D{\Phi}_H(\tau(e),(x_0,\bfm{a}(e))$. 
\item 
The bases of the above three subspaces continuously vary as $(x_0,\bfm{a}(e))$ moves in $\mathscr{M}(\epsilon_2)$. 
\item 
There exist positive constants $C$, $a<1$ such that for all $(x_0,\bfm{a}(e))\in \mathscr{M}(\epsilon_2)$ the following estimates hold. 
\begin{subequations}\label{ineq:nhim_estimates}
\begin{align}
&
\begin{aligned}
&\left\| D{\Phi}_H(\tau(e),(x_0,\bfm{a}(e)))^k u\right\| \leqslant  C(1+|k|)\|u\|\\
    &\qquad \text{for}\ u\in T_{(x_0,\bfm{a}(e))}\mathscr{M}(\epsilon_2), \ k\in\mathbb{Z},
\end{aligned}\\
&
\begin{aligned}
&\left\| D{\Phi}_H(\tau(e),(x_0,\bfm{a}(e)))^k u\right\|\leqslant{}  Ca^k\|u\|\ \\
&\qquad \text{for}\ u\in N^{s}_{(x_0,\bfm{a}(e)))}, \ k\in\mathbb{N},
\end{aligned}\\
&
\begin{aligned}
&\left\| D{\Phi}_H(\tau(e),(x_0,\bfm{a}(e)))^{-k} u\right\|\leqslant{}  Ca^k\|u\|\ \\
    &\qquad \text{for}\ u\in N^{u}_{(x_0,\bfm{a}(e)))}, \ k\in\mathbb{N}.
\end{aligned}
\end{align}
\end{subequations}
\end{enumerate}
\end{theorem}
\begin{proof}
(Step~1) Eigenvalue decomposition of $D{\Phi}_H(\tau(e),(x_0,\bfm{a}(e)))$. \\
Let $|e|\leqslant\epsilon_1$ and $\Gamma(e)$ be the continuation of the trivial solution of (\ref{eqn:hamiltoniansys}) obtained in Proposition~\ref{prop:continuation_for_Hamsys}. Note first that $D{\Phi}_H(\tau(e),(x_0,\bfm{a}(e)))$ is continuous with respect to $(x_0,e)$. 

As in the proof of Proposition~\ref{prop:continuation_for_Hamsys}, using the fact that the monodromy matrix of $\mathrm{Ham}(t)$ has no eigenvalues on the unit circle,  there exists an $\epsilon_3 \in (0,\epsilon_1)$, such that for $|e|\leqslant\epsilon_3$, $D{\Phi}_H(\tau(e),(x_0,\bfm{a}(e)))$ has eigenvalues $1$ with algebraic multiplicity $2$ and two sets of $n-1$ eigenvalues inside and outside the unit circle.  From this eigenvalue distribution, there exist a $2n\times2$ real matrix 
$\Tilde{T}_c^e(0)$ and two $2n\times(n-1)$ matrices $\Tilde{T}_s^e(0)$, $\Tilde{T}_u^e(0)$, which are continuous in $e$, such that $\Tilde{T}(0):=\bmat{\Tilde{T}_c^e(0)&\Tilde{T}_s^e(0)& \Tilde{T}_u^e(0) }$ satisfies 
\[
D{\Tilde{\Phi}}_H(\tau(e),(x_0,\bfm{a}(e)))\Tilde{T}(0)
=\Tilde{T}(0)%
\bmat{\Lambda_c^e& 0& 0 \\ 0&\Lambda_s^e& 0\\ 0& 0& \Lambda_u^e},
\]
where eigenvalues of $\Lambda_c^e$, $\Lambda_s^e$, $\Lambda_u^e$ are $\{1,1\}$, inside 
 and outside of the unit circle, respectively. 
\\[2ex](Step~2) Eigenvalue decomposition of $D\Tilde{\Phi}_H(t,(x_0+\theta,e))$ for $0\leqslant\theta\leqslant 1$. \\
Set $\Tilde{T}^e(t):=D\Tilde{\Phi}_H(t,(x_0,e))\Tilde{T}^e(0)$ and define 
$\Tilde{T}_c^e(\theta)$, $\Tilde{T}_s^e(\theta)$ and $\Tilde{T}_u^e(\theta)$ accordingly as its submatrices. We can prove using (\ref{eqn:shift_monodromy_app}) in Appendix~\ref{sec_app:linearPeriodicSystems} that 
\[
D\Tilde{\Phi}_H(t,(x_0+\theta,e))\Tilde{T}^e(\theta)
=\Tilde{T}^e(\theta)%
\bmat{\Lambda_c^e& 0 & 0 \\ 0 &\Lambda_s^e& 0\\ 0 & 0 & \Lambda_u^e},
\]
for $0\leqslant\theta\leqslant1$. The same decomposition is obtained for $D{\Phi}_H(t,(x_0+\theta,e))$ with a matrix $T^e(\theta)$ made by the reordering linear transformation from $\Tilde{T}^e(\theta)$. 
\\[2ex]
(Step~3) Continuity of the basis. \\
Define 
\begin{gather*}
\Tilde{N}_{(x_0+\theta,e)}^\nu:=\mathrm{Im}\,\Tilde{T}_\nu^e(\theta),\ \nu=\{c,s,u\}
\end{gather*}
and define $N_{(x_0+\theta,e)}^\nu$, $\nu=\{c,s,u\}$ from above by the reordering linear transformation. 
Here, $N_{(x_0+\theta,e)}^c$ corresponds to the tangential part which is given by $T_{(x_0,e)}\mathscr{M}(\epsilon_2)$.
This shows that their bases are continuous in $(x_0,e)$. \\
(Step~4) Take $0<\epsilon_2<\epsilon_3$ and consider $\mathscr{M}(\epsilon_2)$. The estimates (\ref{ineq:nhim_estimates}) are straightforwardly obtained from the decomposition using the eigenvalue properties of $\Lambda_\nu$, $\nu=\{c,s,u\}$ and the boundedness of $\mathscr{M}(\epsilon_2)$. 
\end{proof}

Now, we have the main result in this section. 
\begin{theorem}\label{thm:NHIM}
$\mathscr{M}(\varepsilon_2)$ is a $\rho$-NHIM ($\rho\in\mathbb{N}$) for (\ref{eqn:hamiltoniansys}). 
\end{theorem}
\begin{proof}
The estimates in (\ref{ineq:nhim_contraction_appdx}) in the Appendix~\ref{sctn:nhim_appx} can be verified from (iii) of Theorem~\ref{thm:mainestimates} and one verifies that all the conditions in Definition~\ref{dfn:NHIM_appdx} in the Appendix~\ref{sctn:nhim_appx} for $\mathscr{M}(\epsilon_2)$ to be a $\rho$-\gls{nhim} (for all $\rho\in\mathbb{N}$) for a diffeomorphism induced by the flow of the Hamiltonian system (\ref{eqn:hamiltoniansys}). 
\end{proof}
\begin{remark}
The \gls{nhim} $\mathscr{M}(\epsilon_2)$ has boundaries which are not typically treated in the theory of \gls{nhim} for maps \cite{Hirsch:77:IM}. The treatment of the boundaries will be touched on in the Appendix~\ref{sec_app:boundary}. 
\end{remark}
\section{Existence of the Solution to the Optimal Control Problem}\label{sctn:existence_OC}
Now, based on Theorem~\ref{thm:NHIM}, we are ready to construct a solution for the \gls{hje}~(\ref{eqn:HJE_newCordinate}) from which the optimal control for (\ref{eqn:controlSys_transformed})-(\ref{eqn:cost_new_coordinates}) is obtained. 
Recall that in the theory of \glspl{hje} for optimal point stabilization, the existence of the stable manifold for the extended Hamiltonian system can be analyzed using the theory of algebraic Riccati equations. When the control system is stabilizable in the linear sense, the stable manifold is diffeomorphically projected to the $\bfm{x}$-space via the canonical projection, $(\bfm{x},\bfm{p})\mapsto\bfm{x}$, in $\mathbb{R}^{2n}$ considered as a symplectic manifold. 
After showing that the stable manifold is a Lagrangian submanifold, the stable manifold is represented as a graph of the differential of the solution of the \gls{hje}. 
We refer to \cite{vanderSchaft:91:syscon,vanderSchaft:92:ac} for more details on Lagrangian submanifold and the construction of optimal feedback. 

In the present case, however, this Lagrangian submanifold is not obtained from the stable manifold of the \gls{nhim}, but, from a family of laminations of the \gls{nhim} (see Theorem~\ref{thm:stable_manifold_NHIM_appx} in the Appendix~\ref{sctn:nhim_appx}). It will be shown that a union of specific laminations of $\mathscr{M}(\epsilon_2)$ has a Lagrangian property and satisfies the projectability condition. To this end, we shall look deeper at the structure of the Hamiltonian system (\ref{eqn:hamiltoniansys}) and $N_{(x_0,0)}^s$ (or $\Tilde{T}_s^e(x_0)$) for a point $(x_0,e)$ in $\mathscr{M}(\varepsilon_2)$, listing several useful observations. The following two lemmas show that the bases of vector bundles in Theorem~\ref{thm:mainestimates} have specific properties related to the periodic Riccati differential equation (\ref{eqn:p1Riccati}) (or $\mathrm{Ham}(t)$).

\begin{lemma}
Assume that $M_{\mathrm{Ham}}(1,0)$ has no eigenvalues on the unit circle. Let $\Tilde{U}_\nu(0)$, $\Tilde{V}_\nu(0)$, $\nu\in\{s,u\}$, be real $(n-1)\times(n-1)$ matrices satisfying 
\[
M_{\mathrm{Ham}}(1,0)
\bmat{\Tilde{U}_s(0)& \Tilde{U}_u(0)\\
\Tilde{V}_s(0)& \Tilde{V}_u(0)}
=
\bmat{\Tilde{U}_s(0)& \Tilde{U}_u(0)\\
\Tilde{V}_s(0)& \Tilde{V}_u(0)}
\bmat{\Lambda_s& 0\\ 0 & \Lambda_u},
\]
where $\Lambda_s$ and $\Lambda_u$ have only eigenvalues with inside ($|\lambda| < 1$) and outside ($|\lambda| > 1$) the unit circle, respectively.
Then, there exist unique $(n-1)$-dimensional row vectors $\xi_s$, $\xi_u$ such that $\Tilde{T}^0_s(0)$, $\Tilde{T}^0_u(0)$ defined by
\[
\Tilde{T}^0_\nu(0)=\bmat{
\xi_\nu^\top\\ \bfm{0}_{n-1}^\top\\
\Tilde{U}_\nu(0)\\ \Tilde{V}_\nu(0)
},\ \nu\in\{s,u\}
\]
and $M_{\Tilde{{\mathscr{H}}}}$ (see Lemma~\ref{lemma:form_of_sol_for_H}) satisfy
\[
M_{\Tilde{{\mathscr{H}}}}(1,0)\bmat{\Tilde{T}^0_s(0)& \Tilde{T}^0_u(0)
}=
\bmat{\Tilde{T}^0_s(0)& \Tilde{T}^0_u(0)
}
\bmat{\Lambda_s& 0 \\ 0  & \Lambda_u}.
\]
\end{lemma}
\begin{lemma}\label{fact:MHam}
Consider $\Tilde{T}^e_\nu(\theta)$, $\nu\in\{s,u\}$, constructed in the proof of Theorem~\ref{thm:mainestimates}. At $e=0$, their submatrix from the third to the last row corresponds to 
\[
\bmat{\Tilde{U}_\nu(t)\\\Tilde{V}_\nu(t)
}
=M_{\mathrm{Ham}}(t,0)\bmat{\Tilde{U}_\nu(0)\\\Tilde{V}_\nu(0)
},\ \nu\in\{s,u\}. 
\]
\end{lemma}

\vskip 1ex
The proof of Lemma~\ref{fact:MHam} is a direct computation using Lemma~\ref{lemma:form_of_sol_for_H}. 
The following lemma is shown by Lemma~\ref{fact:MHam} and Proposition~\ref{prop:p-ric_construction_appx} in the Appendix~\ref{sec_app:linearPeriodicSystems} ($X_s(t)$ corresponds to $\Tilde{U}_s(t)$). 
\begin{lemma}\label{fact:stablesubspace_projection}
If $M_{\mathrm{Ham}}(1,0)$ has no eigenvalues on the unit circle and $(A(t),\bar{R}(t))$ is stabilizable, then $\Tilde{U}_s(t)$ is nonsingular for $t\in [0,1]$ and the stabilizing solution to periodic Riccati equation (\ref{eqn:p1Riccati}) exists. 
Moreover, for all $x_0\in [0,1]$ the subspace  $N_{(x_0,0)}^s$ in Theorem~\ref{thm:mainestimates} is isomorphically projectable to the space of $\bfm{x}_a$. 
\end{lemma}

\vskip 1ex
Let $L(x_0) = W^{ss}(x_0,\bfm{0}_{n-1})$ (lamination, see Theorem~\ref{thm:stable_manifold_NHIM_appx} in the Appendix~\ref{sctn:nhim_appx}) for points $(x_0,\bfm{0}_{n-1})$ in $\mathscr{M}(\epsilon_2)$ and let 
\[
L=\bigcup_{(x_0,\bfm{0}_{n-1})\in \mathscr{M}(\epsilon_2)}L(x_0).
\]

\begin{proposition}\label{prop:lagrange_SM_lammination}
Assume that $M_{\mathrm{Ham}}(1,0)$ has no eigenvalues on the unit circle and $(A(t),\bar{R}(t))$ is stabilizable. Then, 
$L$ is a Lagrangian submanifold that is locally diffeomorphically projectable, via the canonical projection, to the $\bfm{x}$-space. 
Moreover, there is an open neighborhood $U\subset \mathbb{R}^n$ of $x_0$-axis and a $C^r$ function $V(\bfm{x})$ that is defined in $U$ and 1-periodic in $x_0$ satisfying HJE (\ref{eqn:HJE_newCordinate}) in $U$, namely, $p_0=\frac{\partial V}{\partial x_0}$ and $\bfm{p}_a=\frac{\partial V}{\partial \bfm{x}_a}^\top$ satisfy (\ref{eqn:HJE_newCordinate}) in $U$.
\end{proposition}
\begin{proof}
First, it will be shown that $L$ is a Lagrangian submanifold, i.e., $\omega=\sum_{i=0}^{2n-1} dx_i\wedge dp_i$ restricted to $L$ vanishes and $L$ has dimension $n$.
From Theorem~\ref{thm:stable_manifold_NHIM_appx},  $L$ is invariant under the Hamiltonian flow of $X_H$ and any solution $\Phi_H(t,q)$ starting in $q= (\bfm x,\bfm p) \in L$ will eventually converge to $\Gamma_0$.
For any $q\in L$ and tangent vectors $Q_1,Q_2\in T_q L$ one obtains that 
\begin{align*}
    \omega(Q_1,Q_2) = \omega(D\Phi_H(t,q) Q_1,D\Phi_H(t,q) Q_2) \text{ for } t\geqslant0.
\end{align*}
As $\Phi_{H}(t, q)$ will converge to $\Gamma_0$, it follows that 
\[D\Phi_H(t,q) Q_i\rightarrow [\ast,\bfm{0}_{2n-1}] \text{ for } t\rightarrow \infty \text{ and } i\in\{1,2\}.\]
Therefore, $\omega(Q_1,Q_2) = 0$ since all components related to $\bfm{p}$ approach zero and $\omega$ vanishes. 
The dimension of $L(x_0)$ is $n-1$ dimensional, consequently $L$ is $n$ dimensional and thus a Lagrangian submanifold.
The local projectability to the $\bfm{x}$-space follows from Lemma~\ref{fact:stablesubspace_projection} and the fact that $L(x_0)$ is tangent to $N_{(x_0,e)}^s$. Thus, by the implicit function theorem, there exists an $\mathbb{R}^n$-valued function $\xi(\bfm{x})$, which is a period-1 function in $x_0$, in some neighborhood $U$ of the $x_0$-axis such that
\[
L=\{ (\bfm{x},\bfm{p})\,|\, \bfm{p}=\xi(\bfm{x}) \}.
\]
The periodicity of $\xi$ is shown as follows. $L(x_0)$ is the set of solutions $\Phi_H(t,q)$ to the Hamiltonian system (\ref{eqn:hamiltoniansys}) converging to $(t+x_0,\bfm{0}_{2n-1})$ as $t\to\infty$ while $L(x_0+1)$ consists of $\Phi_H(t,q)$ converging to $(t+x_0+1,\bfm{0}_{2n-1})$. From the periodicity of (\ref{eqn:hamiltoniansys}) in $x_0$, $L$ is 1-periodic in $x_0$ and so is $\xi$. It is noted that $\xi$ satisfies $\partial \xi_i/\partial x_j=\partial \xi_j/\partial x_i$ for $i,j=0,\cdots,n-1$ from $\omega|_L=0$. 
Using this, as in the standard proof of the (local) Poincar\'e lemma, $V(\bfm{x})$, defined by
\[
V(\bfm{x})=\sum_{i=0}^{n-1}\int_0^{x_i}\xi_i(x_0,x_1,\cdots,x_{i-1},y,\bfm{0}_{n-i-1})\,dy
\]
in a star-shaped subset of $U$, is a period-1 function in $x_0$ and satifies $dV = \xi$, i.e. $V$ is a solution to the HJE (\ref{eqn:HJE_newCordinate}). 
\end{proof}

\vskip 1ex
We are now in the position to state the main theorem of the paper. 
\begin{theorem}
\label{thm:main_theorem}
Suppose the matrices \( A(t), Q(t), \bar{R}(t) \) are given as in \eqref{eqn:cost_new_coordinates} and \eqref{eqn:lin_hamiltoniansys}, and that the pair \( (A(t), \bar{R}(t)) \) is stabilizable while \( (Q(t), A(t)) \) is detectable.
Then, there exits, locally around $\bfm{x}_a=0$, a solution $V(x_0,\bfm{x}_a)$ to the HJE (\ref{eqn:HJE_newCordinate}) such that the solution for Problem~\ref{prblm:main_problem} is locally given, as a feedback control, by
\[
u = -\frac{1}{2}R^{-1}\left(g_0(\bfm{x})^\top \frac{\partial V}{\partial x_0} + g_a(\bfm{x})^\top \frac{\partial V}{\partial \bfm{x}_a}^\top\right). 
\]
\end{theorem}

\section{Application examples}
\label{sec:applications}
In this section, we present two examples. In the first, we illustrate the general procedure and, in particular, demonstrate how to construct the orthogonal coordinates. In the second example, we consider a setting from orbital mechanics, which naturally exhibits transverse coordinates, and show how the proposed control approach can outperform classical linearized control.

\subsection{Energy control for a mass-spring system}
In the following, a mass-spring system is considered, whose dynamics is given by
\begin{align*}
	\dot{z}_1 &= z_2,\\
	\dot{z}_2 &= - z_1 + u,
\end{align*}
where $u$ denotes the system input.
The goal of the control is to stabilize the energy level corresponding to one, i.e., $ \frac{1}{2}(z_1^2+z_2^2)-1 = 0$.
In the following the proposed transversal coordinate system \eqref{eqn:def_trafo_Hale} along the orbit is used. For this, the point transformation 
\[
	x = \begin{bmatrix}x_0\\x_1\end{bmatrix} = \begin{bmatrix}-\arctan(\frac{z_2}{z_1})\\ \frac{1}{2}(z_1^2+z_2^2)-1\end{bmatrix},
\]
is used.
Here, $x_0$ and $x_1$ represent the new coordinates adapted to the periodic orbit of the mass-spring system: $x_0$ denotes the phase angle along the orbit, while $x_1$ measures the distance from the orbit, corresponding to the deviation to the desired energy level.
Please note, that we decided to keep the periodicity of $2\pi$ for the orbit for readability reasons.
In the new coordinates, the dynamics yields
\begin{align*}
	\dot{x}_0&=1 - \cos(x_0)u,\\
	\dot{x}_1&=(2x_1+1)\sin(x_0)u.
\end{align*}
The cost function in form of~\eqref{eqn:cost_new_coordinates} to be minimized is given by
\begin{align*}
    J = \int_0^\infty x_1^2 + u^2 dt,
\end{align*}
where $Q = 2$ and $R=1$.
The HJE \eqref{eqn:HJE_newCordinate} is straightforwardly derived as
\begin{align*}
H\left(\bfm{x},\frac{\partial V}{\partial \bfm{x}}\right)=\frac{\partial V}{\partial \bfm{x}_0}-\frac{1}{4}G\left(\bfm{x},\frac{\partial V}{\partial \bfm{x}}\right)^2+x_1^2=0,
\end{align*}
where $G\left(\bfm{x},\frac{\partial V}{\partial \bfm{x}}\right)=\left(- \frac{\partial V}{\partial \bfm{x}_0} \cos{\left(x_{0} \right)} + \frac{\partial V}{\partial \bfm{x}_1}\left(2 x_{1} + 1\right) \sin{\left(x_{0} \right)}\right)$. 
Moreover, the Hamiltonian flow, i.e., \eqref{eqn:hamiltoniansys}, results in
\begin{align*}
	\dot{x}_0&=1+\frac{1}{2} G(\bfm{x},\bfm{p}) \cos{\left(x_{0} \right)},\\
	\dot{x}_1&=-\frac{1}{2}\left(2 x_{1} + 1\right) G(\bfm{x},\bfm{p}) \sin{\left(x_{0} \right)},\\
	\dot{p}_0&=\!\frac{1}{4}\!\left(2 p_{0} \sin{\left(x_{0} \right)}\!+\!2p_{1}\!\left(2x_{1}\!+\!1\right)\cos{\left(x_{0}\right)}\right)G(\bfm{x},\bfm{p}),\\
	\dot{p}_1&=p_{1} G(\bfm{x},\bfm{p}) \sin{\left(x_{0} \right)} - 2 x_{1}.
\end{align*}
The linear dynamics along the orbit with $x_1=0$, $p_0=0$ and $p_1=0$ is
\begin{align*}
    \mathscr{H}(t)\!=\!
	\begin{bmatrix}0 & 0 & - 0.5 \cos^{2}(t) & 0.5 \sin(t) \cos(t)\\0 & 0 & 0.5 \sin(t)\cos(t)  & - 0.5 \sin^{2}(t)\\0 & 0 & 0 & 0\\0 & -2 & 0 & 0\end{bmatrix}.
 \end{align*}
We read $A(t),Q(t)$ and $\bar{R}(t)$ out of $\mathscr{H}(t)$.
Theorem~\ref{thm:main_theorem} is employed to verify that there exits a solution to the \gls{hje}.
For this, it is required that the pair of $A(t)=0$ and $\bar{R}(t)=\frac{\sin(t)^2}{2}$ is stabilizable.
Moreover, the pair of $Q(t)=2$ and $A(t)=0$ has to be detectable.
While the latter is obvious, the stabilizability is verified by showing that the controllability gramian $W_c(t_0,t_1)$ is invertible for some $t_1>t_0$ in the following.
The controllability gramian for $t_1>t_0$ yields
\begin{align*}
    W_c(t_0,t_1)    & = \int_{t_0}^{t_1} e^{A(t)(t_1-\tau)} \bar{R}(\tau) \bar{R}(\tau)^\top {e^{A(t)(t_1-\tau)}}^\top d\tau,\\
                    & = \frac{1}{4}\int_{t_0}^{t_1} \sin^4(\tau)  d\tau>0,
\end{align*}
and it follows that $(A(t),\bar{R}(t))$ is controllable (and therefore stabilizable).
Consequently, locally near $\Gamma_0=\{(\bfm{x},\bfm{p})\in T \mathbb{R}^2|x_1=0, \bfm{p}=0\}$ there exists a solution of the HJE \eqref{eqn:HJE_newCordinate} which guarantees the existence of a feedback law near $\Gamma_0$.
The fundamential matrix solution $D\Phi_H(t,\bfm{q})$ for $\bfm{q}=(\bfm x, \bfm p)$ was computed numerically for $t = 2 \pi$ (after one full period of the orbit) and unit matrix as initial condition and yields
\begin{align*}
D\Phi_H(2\pi,\bfm{q}) \approx \begin{bmatrix}
1			&15.3	& 0.3  	& -4.3\\
0  			& 38.7	& 4.3	& -11.2\\
0  			& 0			& 1			& 0		\\
0 			& -133.2	&  -15.3	& 38.7	
\end{bmatrix}.
\end{align*}
As predicted by Proposition~\ref{prop:pRic_eigenvalues}, the system exhibits an eigenvalue equal to one with algebraic multiplicity two, along with a pair of reciprocal eigenvalues: \(7.7336 \times 10^1\) and \(1.2931 \times 10^{-2}\).

The optimal control problem was solved using numerical optimization via a direct collocation method with 1000 knot points over a time horizon of 10 seconds. The resulting nonlinear program was solved using the Interior Point OPTimizer (IPOPT), a software library for large-scale nonlinear optimization.
Figure \ref{fig:spring_mass_system} shows the resulting trajectories in the $(z_1,z_2)$-plane, which approximate the projection of the Hamiltonian flow along $L$ onto the canonical coordinates. The corresponding control inputs for representative initial conditions are given in Fig.~\ref{fig:mass_spring_control_inputs}. For comparison, we also computed the time-varying \gls{pre}-based linearized controller (dashed curves); the associated performance values are listed in Table~\ref{tab:costs_mass_spring}. For initial conditions far from the orbit, the nonlinear feedback yields markedly superior performance relative to the linearized controller.
\begin{figure}
    \centering
    \input{graphics/mass_spring_system_state}
    \caption{Solid trajectories correspond to the Hamiltonian flow along the Lagrangian submanifold $L$ projected onto $z_1$-$z_2$ plane. The dashed lines correspond to the \gls{pre}-based linear feedback.
}
    \label{fig:spring_mass_system}
    \end{figure}
\begin{figure}
    \input{graphics/mass_spring_system_input}
	\caption{Control inputs for selected initial states. Solid lines denote nonlinear feedback; dashed lines denote linear \gls{pre}-based feedback.}
	\label{fig:mass_spring_control_inputs}
\end{figure}
\begin{table}
\centering
\begin{tabular}{c c c c c}
\toprule
$x_1(0)$ & $x_2(0)$ & $J_{\text{nonlinear}}$ & $J_{\text{linear}}$ & $\frac{J_{\text{nonlinear}}}{J_{\text{linear}}}$ \\
\midrule
0.000 & 10.000 & 329.4873 & 337.6515 & 0.9758 \\
5.000 & 8.660 & 664.5139 & 789.1495 & 0.8421 \\
8.660 & 5.000 & 1681.1021 & 3036.7890 & 0.5536 \\
10.000 & 0.000 & 2130.0929 & 4368.4101 & 0.4876 \\
8.660 & -5.000 & 1107.9184 & 1629.4416 & 0.6799 \\
5.000 & -8.660 & 390.2213 & 392.1030 & 0.9952 \\
-0.000 & 0.500 & 1.0717 & 1.0744 & 0.9975 \\
0.250 & 0.433 & 1.2785 & 1.2814 & 0.9978 \\
0.433 & 0.250 & 1.6965 & 1.7638 & 0.9618 \\
0.500 & 0.000 & 1.2963 & 1.3396 & 0.9676 \\
0.433 & -0.250 & 1.0955 & 1.0983 & 0.9975 \\
0.250 & -0.433 & 1.0202 & 1.0216 & 0.9986 \\
\bottomrule
\end{tabular}
\caption{Initial conditions, corresponding costs, and cost ratios between the nonlinear and linear control strategies for the mass spring energy scenario.}
\label{tab:costs_mass_spring}
\end{table}

\subsection{Satellite orbit transfer}
In this subsection, optimal feedback control is applied to an orbital mechanics setting.
The dynamics is that of a massless body moving in a central gravitational force field subject also to drag and a radial modulated force, which was considered in \cite{Owis2009}.
The equations of motion can be stated as
\begin{align*}
    \dot{z}_1 &= z_3\\
    \dot{z}_2 &= \frac{1-\gamma z_2}{z_1^2}\\
    \dot{z}_3 &= \frac{(1-\gamma z_2)^2}{z_1^3}- \frac{\gamma z_3 + 1}{z_1^2} + u,
\end{align*}
where $z_1$ and $z_2$ denote the radial and angular components of the polar coordinates, respectively.
Furthermore, $z_3$ is the radial velocity, $\gamma$ denotes the drag coefficient and $u$ is the control input.
As in~\cite{Owis2009}, the coordinate representing the angular velocity, i.e., $\dot{z}_2$, was eliminated by exploiting
the symmetry introduced by the conserved quantity
\[
h = z_1^2 \dot{z}_2 + \gamma z_2.
\]
The goal of the control is to stabilize the desired orbit determined by the radius $z_1=1$.
Analyzing the conserved quatity, it becomes obvious that a system with drag, i.e. $\gamma>0$, inevitably loses angular momentum $z_1^2 \dot{z}_2$ and the tangential velocity $\dot{z}_2$ approaches zero. 
Therefore, the periodic orbit reduces to a point. 
Since we want to study an actual periodic orbit, a drag-free system, i.e. $\gamma$ equals to zero, is assumed, which is a reasonable assumption for typical use cases. 
Here, the required input along the desired trajectory, given by $z_1=1$, $z_2(t)=t+\text{const}$ and $z_3(t)=0$ is zero and the invariance condition for the unforced system is satisfied.


As a next step, a change of coordinates is performed to obtain the system dynamics in transverse coordinates defined by
\begin{align*}
    \begin{bmatrix}
        x_0\\
        x_1\\
        x_2
    \end{bmatrix} = \begin{bmatrix}
        z_2\\
        z_1-1\\
        z_3
    \end{bmatrix}.
\end{align*}
Here, $x_0$ denotes the phase of the orbit, represented by the angular component in polar coordinates. The variable $x_1$ captures the radial deviation from the orbit, while $x_2$ measures the deviation in radial velocity—$x_2$ is zero along the desired periodic orbit. We have chosen to retain the orbit’s natural periodicity of $2\pi$ rather than normalizing it to 1, in order to enhance readability.
The system dynamics results in
\begin{align*}
    \dot{x}_0 &= \frac{1}{(x_1+1)^2} = 1 +f_0(x_1)\\
    \dot{x}_1 &= x_2\\
    \dot{x}_2 &= \frac{1}{(x_1+1)^3} - \frac{1}{(x_1+1)^2} + u = -x_1 +f_2(x_1) + u,
\end{align*}
where $f_0(x_1) = - 2 x_1 + 3 x_1^2 + \mathcal{O}(|x_1|^3)$ and $f_2(x_1) = 3 x_1^2 + \mathcal{O}(|x_1|^3)$.
In contrast to the finite time horizon cost of \cite{Owis2009}, here an infinite horizon cost function 
\begin{align*}
    \int_0^\infty \frac{1}{20} x_1^2+u^2 dt
\end{align*}
is subject to minimization, such that the satellite follows an orbit at a constant altitude of $z_1=1$.

The \gls{hje} in the new coordinates is
\begin{multline*}
    H\left(\bfm{x},\frac{\partial V}{\partial \bfm{x}}\right)= \frac{\partial V}{\partial \bfm{x_0}}+\frac{\partial V}{\partial \bfm{x}_0}f_0(x_1)-\frac{1}{4}\left(\frac{\partial V}{\partial \bfm{x}_2}\right)^2\\-\frac{\partial V}{\partial \bfm{x}_1} x_2+\frac{\partial V}{\partial \bfm{x}_2} (-x_1 + f_2(x_1))+\frac{1}{20} x_1^2=0.
\end{multline*}
Please note that the Hamiltonian $H$ is independent of $x_0$ and thus the system has another symmetry with conserved quantity $p_0=\frac{\partial V}{\partial x_0}$.
The Hamiltonian flow is given by
\begin{align*}
    \dot{x}_0 &= 1 + f_0(x_1)\\
    \dot{x}_1 &= x_2\\
    \dot{x}_2 &= -\frac{p_2}{2}-x_1+f_2(x_1)\\
    \dot{p}_0 &= 0 \\
    \dot{p}_1 &= -p_0 \frac{\partial f_0}{\partial x_1}-p_2\frac{\partial f_2}{\partial x_1} + p_2- \frac{1}{10}x_1\\
    \dot{p}_2 &= -p_1.
\end{align*}
The linear dynamics along the orbit with $x_1=0$, $x_2=0$, $p_0=0$, $p_1=0$ and $p_2=0$ is
\begin{align*}
	\mathscr{H}(t)\!=\!
	\begin{bmatrix}0 & -2 & 0 & 0 & 0 & 0\\0 & 0 & 1 & 0 & 0 & 0 \\ 0 & -1 & 0 & 0 & 0 & -\frac{1}{2} \\0 & 0 & 0 & 0 & 0 & 0\\ 0 & -\frac{1}{10}& 0 & 2 & 0 & 1 \\ 0 & 0 & 0 & 0 & -1 & 0\end{bmatrix}.
\end{align*}
The submatrix $\text{Ham}(t)$ is characterized by
\begin{align*}
    \bar{A} = \begin{bmatrix}
        0 & 1\\ -1 & 0
    \end{bmatrix}, && \bar{R}=\begin{bmatrix}
        0 & 0\\ 0 & \frac{1}{2}
    \end{bmatrix} ,&& Q = \begin{bmatrix}
        \frac{1}{10} & 0\\ 0 & 0
    \end{bmatrix}.
\end{align*}
The solution to the periodic Riccati equation reduces to its algebraic counterpart and is given by
\begin{align*}
    P \approx \begin{bmatrix}
    0.36 &   0.09 \\
    0.09 &    0.30
    \end{bmatrix}.
\end{align*}
Therefore, locally a stabilizing control law exists and the solution of the \gls{hje} may be approximated by
 \begin{align*}
    V(\bfm x) = \frac{1}{2} \bfm x_a^\top P \bfm x_a + \mathcal{O}(|\bfm{x}_a|^3).
\end{align*}
Please note, that the Riccati equation as well as the value function is independent of $x_0$ due to the symmetry with conserved quantity $p_0$.
This symmetry also allows to approximate the proposed control law using the stable manifold method \cite{Sakamoto:08:ieee-ac,Sakamoto:13:automatica}.
To show the benefits of the proposed control law, we compare it to the corresponding linearized quadratic optimal controller given by $u_\text{lin.} = -\frac{1}{2}\bar{R}^{-1} g(x)^\top P \bfm x_a$, which is the linear counterpart of Theorem~\ref{thm:main_theorem}.
Both controllers exhibit similar performance near the periodic orbit. However, farther from the orbit, the proposed nonlinear control law significantly outperforms the linearized controller in terms of cost. See Fig.~\ref{fig:satellite_control_state} and Fig.~\ref{fig:satellite_control_inputs} for results corresponding to a representative initial condition, and TABLE~\ref{tab:cost_ratios} for a comparison across various initial conditions.

\begin{figure}
	\input{graphics/satellite_orbit_state}
	\caption{Plot of the state $x$ for selected initial conditions.}
	\label{fig:satellite_control_state}
\end{figure}

\begin{figure}
%
%
\definecolor{mycolor1}{rgb}{0.00000,0.44700,0.74100}%
\definecolor{mycolor2}{rgb}{0.85000,0.32500,0.09800}%
\begin{tikzpicture}

\begin{axis}[%
width=0.85\linewidth,
height=1.5cm,
at={(0.758in,1.511in)},
scale only axis,
xmin=0,
xmax=17,
xlabel style={font=\color{white!15!black}},
xlabel={Time in seconds},
ymin=-6,
ymax=1,
ylabel style={font=\color{white!15!black}},
ylabel={Control},
axis background/.style={fill=white},
xmajorgrids,
ymajorgrids,
legend style={at={(0.97,0.03)}, anchor=south east, legend cell align=left, align=left, draw=white!15!black}
]
\addplot [color=mycolor1, dashed, line width=1.5pt]
  table[row sep=crcr]{%
0	-2.73242744769659\\
0.191922112023704	-2.54870027056637\\
0.405168903161155	-2.35117904762508\\
0.51179229872988	-2.25489937563307\\
0.955167956331692	-1.87295158443806\\
1.39854361393351	-1.52295211273681\\
1.84191927153532	-1.20726448110541\\
2.28529492913713	-0.92733538163781\\
2.69967278407589	-0.698383291030403\\
3.11405063901466	-0.499992058666813\\
3.52842849395342	-0.330573118251177\\
3.94280634889219	-0.18818848625995\\
4.34747697123208	-0.0730117972469557\\
4.75214759357198	0.021071023814784\\
5.15681821591188	0.0966296640643378\\
5.56148883825178	0.156295382822943\\
6.09986077795074	0.215358703594614\\
6.6382327176497	0.257080363315424\\
7.17660465734867	0.286711713630591\\
7.90948716682882	0.312931292041633\\
8.29850830639121	0.318660345218778\\
8.4930188761724	0.316842326519588\\
8.6875294459536	0.309366696940227\\
8.88204001573479	0.292963201531308\\
9.07655058551599	0.263849327947714\\
9.27106115529718	0.217881458168492\\
9.37712992030538	0.184575899783052\\
9.48319868531357	0.146215785727264\\
9.69533621532996	0.0602766188439148\\
9.80140498033816	0.0169685789606504\\
9.90747374534635	-0.023310532320501\\
10.0135425103545	-0.0587217717017907\\
10.1196112753627	-0.0882185904642476\\
10.2110546862568	-0.108648330796697\\
10.3024980971508	-0.124580714114085\\
10.3939415080448	-0.136341187610864\\
10.4853849189388	-0.144401945962983\\
10.6029106627039	-0.150178375001413\\
10.720436406469	-0.151698024118936\\
10.8379621502341	-0.149858718702678\\
10.9554878939992	-0.14544716061781\\
11.1182394123593	-0.136251708683499\\
11.4437424490794	-0.111270226547507\\
12.4835352017215	-0.0209026383701172\\
12.9638996614407	0.0133204498459136\\
13.2040818913002	0.027217384087681\\
13.4442641211598	0.0385786399999368\\
13.6902367082011	0.0472875551580181\\
13.9362092952424	0.0528311937426977\\
14.1821818822838	0.0550685542194991\\
14.4281544693251	0.0538972170147467\\
14.8262979235456	0.045137511207006\\
15.2244413777661	0.0300246388201089\\
16.0581278740006	-0.00555427039849477\\
16.502639620228	-0.0177484940955956\\
16.855695908133	-0.0222855141235172\\
17.0322240520855	-0.0229516934270428\\
};
\addlegendentry{$u(t)$ linear control}

\addplot [color=mycolor2, line width=1.5pt]
  table[row sep=crcr]{%
0	-5.53692920233027\\
0.139214951455706	-5.01620285634494\\
0.269392528237795	-4.54752716976685\\
0.395753354677591	-4.11117489357064\\
0.518010622150474	-3.70777990116284\\
0.625078384831614	-3.37053180622596\\
0.72928523216482	-3.05722580818088\\
0.835194189908258	-2.75428589760275\\
0.932651086497309	-2.48953858804743\\
1.0323588855166	-2.232699483839\\
1.12316782728576	-2.0111485326272\\
1.21631984990186	-1.79608772745789\\
1.31191528271644	-1.58813429094997\\
1.39830816067751	-1.41117287782603\\
1.48718798243417	-1.23981772323705\\
1.5792332853499	-1.07357917009493\\
1.66027046826742	-0.936463788396761\\
1.74245401962952	-0.806032349311383\\
1.82587717560235	-0.682282614390527\\
1.91058403974282	-0.565284664992792\\
1.99616822804349	-0.455648382561542\\
2.08248964660874	-0.353490430566456\\
2.16883384027792	-0.259446392107041\\
2.25426934128532	-0.174080554866244\\
2.33883866175646	-0.096788338970466\\
2.42245982128106	-0.0270953473613424\\
2.50527851086414	0.0356429420001909\\
2.58731183117655	0.0919254729606607\\
2.6685714017207	0.142225201951167\\
2.74906395130312	0.186991196720392\\
2.84186071994215	0.232725867439999\\
2.93366906266666	0.272190451289472\\
3.02436574719152	0.305939924120491\\
3.1266583489293	0.338258668738018\\
3.2274308954229	0.364667581368259\\
3.3267492580838	0.385937377228242\\
3.43650482396957	0.404539861305913\\
3.55837528381403	0.419911113138646\\
3.6810937800998	0.430570593872524\\
3.81724416792986	0.437654639629748\\
3.96738145396981	0.440683398105829\\
4.13185778131159	0.439100055359269\\
4.28625493600196	0.432939191954414\\
4.41701728510077	0.423063047812416\\
4.52072899669176	0.410561776528226\\
4.60742153044876	0.395267634604295\\
4.67887759599754	0.377929624542578\\
4.74329965072035	0.357404008296694\\
4.80268740163308	0.333304116883081\\
4.85855648733013	0.305194690522619\\
4.91233403794872	0.272457835044683\\
4.96530470301184	0.234330636307391\\
5.02160613047795	0.187350886556573\\
5.08368875416938	0.128545234128087\\
5.16256078809808	0.0462284068200312\\
5.31898429851449	-0.118656069150077\\
5.38275589908774	-0.177401624675323\\
5.43820141999911	-0.221854710856192\\
5.49010542181566	-0.25738832506115\\
5.54359083477037	-0.287858768480163\\
5.5957592724347	-0.311874513377195\\
5.65152028617617	-0.331853691777425\\
5.7073455926543	-0.346624908078823\\
5.76883173401098	-0.3576591595574\\
5.8317372117779	-0.364143334157607\\
5.90335687243477	-0.366649945257883\\
5.98002294245629	-0.364702751768554\\
6.07116787323158	-0.357552938236488\\
6.17779049406485	-0.344242078251465\\
6.31173119644462	-0.322320226350119\\
6.49083796459302	-0.28757114020901\\
6.83753808654762	-0.213851456408221\\
7.18193067371685	-0.142717375076156\\
7.43429015882652	-0.0956076566843898\\
7.66137354378098	-0.0579211776119237\\
7.89420314472227	-0.0242848824899866\\
8.11453626621006	0.00280693013510813\\
8.35439464121986	0.0271581679216233\\
8.59375673713534	0.046322551876802\\
8.84459242096462	0.0611676416525491\\
9.10022373838744	0.0710402598400535\\
9.35544393621406	0.0757933918473057\\
9.60813412557131	0.0755872793095733\\
9.87211178804216	0.0702698550235681\\
10.1376573785216	0.0599645357977607\\
10.4382655311376	0.0432640691230652\\
10.9491739231557	0.00872863989483363\\
11.3332047940633	-0.0149216765903049\\
11.6228819975587	-0.0279564651243867\\
11.8957909142505	-0.035489937441568\\
12.1914824561084	-0.0385722349360442\\
12.5089743841455	-0.0369349590334949\\
12.9117611796685	-0.0297077527403573\\
13.6914533801347	-0.0093509588188887\\
14.3520214448857	0.00541526370022893\\
14.8732475383217	0.0121630642972157\\
15.3943539556912	0.0140087928077399\\
15.9953558614282	0.0110468133750565\\
17.0055741821751	0.000772989758431208\\
};
\addlegendentry{$u(t)$ nonlinear control}

\end{axis}
\end{tikzpicture}%
	\caption{Plot of the control input for selected initial conditions.}
	\label{fig:satellite_control_inputs}
\end{figure}

\begin{table}
\centering
\begin{tabular}{c c c c c}
\toprule
$x_1(0)$ & $x_2(0)$ & $J_{\text{nonlinear}}$ & $J_{\text{linear}}$ & $\frac{J_{\text{nonlinear}}}{J_{\text{linear}}}$ \\
\midrule
4.000 & -3.296 & 1.5908 & 1.6966 & 0.9377 \\
3.524 & -4.000 & 1.9163 & 2.3341 & 0.8210 \\
4.000 & -2.400 & 1.3844 & 1.4646 & 0.9452 \\
1.172 & -4.000 & 2.2737 & 3.7794 & 0.6016 \\
2.337 & 4.000 & 8.3166 & 15.4993 & 0.5366 \\
0.544 & -4.000 & 2.5666 & 4.5357 & 0.5659 \\
1.698 & 4.000 & 7.0323 & 13.2758 & 0.5297 \\
-0.332 & 4.000 & 3.9678 & 7.6939 & 0.5157 \\
\bottomrule
\end{tabular}
\caption{Summary of initial conditions, corresponding costs, and cost ratios between the nonlinear and linear control strategies for the satellite orbit transfer scenario.}
\label{tab:cost_ratios}
\end{table}

\section{Conclusion and outlook}
\label{sec:conclusion}
We have shown that an optimal orbital stabilizing controller exists locally provided that stabilizability and detectability conditions along the orbit are satisfied. This has been done by showing the existence of a local solution to an \gls{hje} using the continuation technique, the theories of \gls{nhim} and periodic differential Riccati equations, and symplectic geometry (Hamiltonian mechanics). This work extends a well-known result in the nonlinear point stabilization problem to orbital stabilization and can be considered, especially, as the generalization of the work in \cite{vanderSchaft:91:syscon,vanderSchaft:92:ac} in the sense that the framework employed in the present paper purely generalizes that in \cite{vanderSchaft:91:syscon,vanderSchaft:92:ac}. 
The computational aspect, however, has not been fully addressed in the paper. There are several challenges. The first is to obtain the system description (\ref{eqn:controlSys_transformed}) in moving orthogonal coordinates. The second is to compute a solution in a periodic Riccati differential equation. The work in \cite{Gusev2010} may be useful in this respect. The third is to compute the union of lamination $L$ in Proposition~\ref{prop:lagrange_SM_lammination}, which plays the same role as the stable manifold of the Hamiltonian extension in \cite{Sakamoto:08:ieee-ac} and \cite{Sakamoto:13:automatica}. 

Future work will aim to extend existing numerical methods for point stabilization~\cite{Sakamoto:08:ieee-ac,Sakamoto:13:automatica} to the orbital stabilization setting. It would be particularly interesting to investigate how nonlinearities in underactuated mechanical systems can be exploited and how they influence optimal orbital stabilization behavior.
Another interesting direction for future research is the extension of the proposed framework to hybrid dynamical systems, which naturally arise in robotic locomotion involving impacts, e.g., \cite{Westervelt2003}. Under suitable conditions, it may be possible to show that the associated stable manifold persists by approximating the hybrid dynamics in a continuous form. 

The last topic in the outlook is the possibility of another way of proving the existence of optimal control reported in \cite{Sakamoto:22:automatica} using a nonlinear functional analysis technique. The conditions obtained there are given in terms of exponential stabilizability, detectability, and nonlinear growth conditions, thus potentially suitable for orbital stabilization using the prior works in the orbital stabilization mentioned in the Introduction (\S~\ref{sctn:introduciton}). 

\section*{References}
\bibliographystyle{IEEEtran}

%
\appendices
\section{Theory of linear periodic systems and differential Riccati equations}
\label{sec_app:linearPeriodicSystems}
\setcounter{equation}{0}
\renewcommand{\theequation}{\Roman{section}.\arabic{equation}}
\setcounter{theorem}{0}
\renewcommand{\thetheorem}{\Roman{section}.\arabic{theorem}}
\setcounter{proposition}{0}
\renewcommand{\theproposition}{\Roman{section}.\arabic{proposition}}
\setcounter{remark}{0}
\renewcommand{\theremark}{\Roman{section}.\arabic{remark}}
\setcounter{factnew}{0}
\renewcommand{\thefactnew}{\Roman{section}.\arabic{factnew}}
\setcounter{notation}{0}
\renewcommand{\thenotation}{\Roman{section}.\arabic{notation}}
\setcounter{definition}{0}
\renewcommand{\thedefinition}{\Roman{section}.\arabic{definition}}
%

In this section, we review basic facts on linear periodic system theory such as periodic differential Riccati equations. All the materials presented in this Appendix~\ref{sec_app:linearPeriodicSystems} can be found in \cite{Kano1979:ijc,Kano1985:ieeetac,Ichikawa:01:LTVSSS}. A matrix-valued function $U(t)$, $t\in\mathbb{R}$, is said to be $t_p$-periodic if $U(t)=U(t+t_p)$ for $t\in\mathbb{R}$. 

Let $A(t)$ be a $t_p$-periodic real matrix of $n\times n$ dimensions. Let also  $M_A(t,s)$ be the state transition matrix for the differential equation $\dot{x}(t)=A(t)x(t)$, namely, 
\[
\frac{\partial}{\partial t}M_A(t,s)=A(t)M_A(t,s),
\quad M_A(t,t)=I.
\]
The Floquet theory (see, e.g., \cite[p.117]{Hale:73:ODE}) says that $M_A(t,s)$ satisfies $t_p$-periodicity $M_A(t+t_p,s+t_p)=M_A(t,s)$ for $t,s\in \mathbb{R}$. $M_A(t_p,0)$ is called the monodromy matrix and plays a key role in analyzing linear periodic differential equations. For instance, the following eigenvalue decomposition will be used on many occasions
\[
M_A(t_p,0)\bmat{X_1&X_2\\ Y_1&Y_2}=\bmat{X_1&X_2\\ Y_1&Y_2}\bmat{\Lambda_1&0\\0&\Lambda_2}, 
\]
where all the matrices above are real and $\Lambda_1$, $\Lambda_2$ do not have common eigenvalues. It can also be shown that, $X_j(t)$, $Y_j(t)$ ($j=1,2$) defined by
\[
\bmat{X_1(t)&X_2(t)\\ Y_1(t)&Y_2(t)}=M(t,0)\bmat{X_1&X_2\\ Y_1&Y_2}
\]
satisfy 
\begin{subequations}\label{eqn:shift_monodromy_app}
\begin{gather}
\frac{d}{dt}\bmat{X_1(t)&X_2(t)\\ Y_1(t)&Y_2(t)}=A(t)\bmat{X_1(t)&X_2(t)\\ Y_1(t)&Y_2(t)}\label{eqn:eig_deomp_differential_eq_app}\\
M_A(t+t_p,t)\bmat{X_1(t)&X_2(t)\\ Y_1(t)&Y_2(t)}
    =\bmat{X_1(t)&X_2(t)\\ Y_1(t)&Y_2(t)}\bmat{\Lambda_1&0\\0&\Lambda_2}.
\label{eqn:shifted_eig_docmp_app}
\end{gather}
\end{subequations}
%
%

\begin{definition}
A $t_p$-periodic square matrix $A(t)$ is said to be asymptotically stable if the corresponding linear differential equation $\dot{x}=A(t)x$ is asymptotically stable at the origin. 
\end{definition}

\vskip 1ex
Now, let $B(t)$, $C(t)$ be $t_p$-periodic real matrices of $n\times m$ and $r\times n$ dimensions, respectively and consider a $t_p$-periodic linear control system
\begin{equation}
\dot{x}=A(t)x+B(t)u, \ y=C(t)x, \label{eqn:periodic_linearsys}
\end{equation}
where $u(t)\in\mathbb{R}^m$ is the control input and $y(t)\in\mathbb{R}^r$ is the output. The stabilizability and detectability properties, which play an important role in linear time-invariant control systems, are defined for (\ref{eqn:periodic_linearsys}) as follows. 
\begin{definition}
\begin{enumerate}[(i)]
    \item An eigenvalue $\lambda$ of $M_A(t_p,0)$ is said to be $(A(t),B(t))$-controllable if $M_A(t_p,0)^\top \xi=\lambda\xi$ and $B(t)^\top M_A(t,0)^{-\top}\xi =0$ for $t\in [0,t_p]$ imply $\xi=0$. 
    \item An eigenvalue $\lambda$ of $M_A(t_p,0)$ is said to be $(C(t),A(t))$-observable if $M_A(t_p,0) \xi=\lambda\xi$ and $C(t)\Phi_A(t,0)\xi =0$ for $t\in [0,t_p]$ imply $\xi=0$.
    \item The pair $(A(t),B(t))$ is said to be stabilizable if all eigenvalues $\lambda$ of $M_A(t_p,0)$ with $|\lambda|\geqslant 1$ are $(A(t),B(t))$-controllable.
    \item The pair $(C(t),A(t))$ is said to be detectable if all eigenvalues $\lambda$ of $M_A(t_p,0)$ with $|\lambda|\geqslant 1$ are $(C(t),A(t))$-observable. 
    \item The pair $(A(t),B(t))$ is said to be controllable if all eigenvalues of $M_A(t_p,0)$ are $(A(t),B(t))$-controllable. 
    \item The pair $(C(t),A(t))$ is said to be detectable if all eigenvalues of $M_A(t_p,0)$ are $(C(t),A(t))$-observable.
\end{enumerate} 
\end{definition}

\vskip 1ex
It is shown that $(A(t),B(t))$ is stabilizable if and only if there exists a continuous $t_p$-periodic $m\times n$ matrix $K(t)$ such that $A(t)+B(t)K(t)$ is asymptotically stable and that $(C(t), A(t))$ is detectable if and only if there exists a continuous $T$-periodic $n\times r$ matrix $G(t)$ such that $A(t)+G(t)C(t))$ is asymptotically stable.

\vskip 1ex
The periodic Riccati equation, which plays the central role in optimal control for (\ref{eqn:periodic_linearsys}), takes the following form \begin{equation}
-\dot{P}(t)=P(t)A(t)+A(t)^\top P(t)-P(t)R(t)P(t)+Q(t),\label{eqn:periodicRic_appx}
\end{equation}
where $R(t)$, $Q(t)$ are $t_p$-periodic positive semi-definite matrices of $n\times n$ dimension. 
\begin{theorem}\label{thm:periodic_riccati_existence_appx}
A necessary and sufficient condition for (\ref{eqn:periodicRic_appx}) to have a solution $P(t)$ with $A(t)-R(t)P(t)$ being asymptotically stable (stabilizing solution) is that $(A(t),R(t))$ is stabilizable and all eigenvalues of $M_{A}(t_p,0)$ on the unit circle, if they exist, are $(Q(t),A(t))$-detectable. Under this condition, $P(t)$ is a $t_p$-periodic positive semi-definite matrix, and no other solution exists with the closed-loop stability. 
\end{theorem}

\vskip 1ex
The construction of the stabilizing solution is done as follows. Let $M_{\mathrm{Ham}}(t,s)$ denote the state transition matrix of the Hamiltonian matrix
\[
\mathrm{Ham}=\bmat{A(t)& -R(t)\\ -Q(t)& -A(t)^\top}.
\]
If the monodromy matrix $M_{\mathrm{Ham}}(t_p,0)$ has no eigenvalues on the unit circle, we have a decomposition
\[
M_{\mathrm{Ham}}(t_p,0)\bmat{X_s(0)& X_u(0)\\ Y_s(0)& Y_u(0)}=
\bmat{X_s(0)& X_u(0)\\ Y_s(0)& Y_u(0)}\bmat{\Lambda_s&0 \\0&\Lambda_u}
\]
with appropriate real matrices, where the eigenvalues of $\Lambda_s$ are inside the unit circle whereas those of $\Lambda_u$ are outside. Set
\[
\bmat{X_s(t)& X_u(t)\\ Y_s(t)& Y_u(t)}=
M_{\mathrm{Ham}}(t,0)\bmat{X_s(0)& X_u(0)\\ Y_s(0)& Y_u(0)}. 
\]
\begin{proposition}\label{prop:p-ric_construction_appx}
If $M_{\mathrm{Ham}}(t_p,0)$ has no eigenvalues on the unit circle and $(A(t),R(t))$ is stabilizable, then $\det(X_s(t))\neq0$ for $t\in[0,t_p]$ and the stabilizing solution of (\ref{eqn:periodicRic_appx}) is given by $P(t)=Y_s(t)X_s(t)^{-1}$. 
\end{proposition}

\vspace{1ex} 
The $t_p$-periodicity of $P(t)$ is shown using the property of $M_{\mathrm{Ham}}(t,s)$ or using the technique in \cite[Theorem 2.2]{Ichikawa:01:LTVSSS}. 
%
\section{Normally hyperbolic invariant manifolds}\label{sctn:nhim_appx}
\setcounter{equation}{0}
Let $F:\mathbb{R}^N\to\mathbb{R}^N$ ($N\geqslant3$) be a $C^r$ ($r\geqslant1$)  diffeomorphism and $M\subset\mathbb{R}^N$ be a compact $C^r$ submanifold which is invariant under $F$; $F(M)=M$. 
\begin{definition}\label{dfn:NHIM_appdx}
The invariant manifold $M$ is called {\em $r$-normally hyperbolic invariant manifold} (NHIM) for $F$ if the following hold. 
\begin{enumerate}[(i)]
    \item 
There exists a continuous $DF$-invariant splitting 
\begin{equation}
\left. T\mathbb{R}^N\right|_M=TM\oplus N_M^s\oplus N_M^u\label{eqn:split_appdx}
\end{equation}
of the tangent bundle $T\mathbb{R}^N$ over $M$ such that 
\begin{gather*}
DF(x)(T_xM)=T_{F(x)}M,\ DF(x)(N_x^s)=N_{F(x)}^s, \\
\ DF(x)(N_x^u)=N_{F(x)}^u
\end{gather*}
hold for $x\in M$. Here, continuity of the split means that as $x$ varies in $M$ one can find continuously varying bases in $N_x^s$ and $N_x^u$. 
\item
there exist positive constants $C$ and $0<a<1$ such that 
\begin{subequations}\label{ineq:nhim_contraction_appdx}
\begin{align}
&
\left\| \left.(DF)^k\right|_{N_x^s} \right\|\cdot
    \left\| \left( \left. (DF)^k \right|_{T_xM}  \right)^{-1} \right\|^\rho < C a^k\\
&
\left\| \left( \left.(DF)^k\right|_{N_x^u} \right)^{-1} \right\|\cdot
    \left\|  \left. (DF)^k \right|_{T_xM}   \right\|^\rho < C\ a^k
\end{align}
\end{subequations}
for $0\leqslant\rho\leqslant r$, $x\in M$ and $k\in\mathbb{N}$, where $\|\cdot\|$ stands for the induced norm for a linear map. 
\end{enumerate}
\end{definition}

\vskip 1ex
We now state the stable manifold theorem for $M$ (see \cite{Hirsch:77:IM} also see \cite{Wiggins:94:NHIMDS} for flows). 
\begin{theorem}\label{thm:stable_manifold_NHIM_appx}
If $M$ is an NHIM for $F$, then, local stable manifold $W^s(M)$ and local unstable manifold $W^u(M)$ exist, which are $C^r$ submanifolds of $\mathbb{R}^N$. $W_{\mathrm{loc}}^s(M)$ and $W_{\mathrm{loc}}^u(M)$ are tangent to $TM\oplus N^s$ and $TM\oplus N^u$, respectively, at each point of $M$. Moreover, there exist two $F$-invariant laminations $W_{\mathrm{loc}}^{ss}(x)$ and $W_{\mathrm{loc}}^{uu}(x)$ ($x\in M$), which are leaves of  $W_{\mathrm{loc}}^s(M)$ and $W_{\mathrm{loc}}^u(M)$, respectively, and defined as 
\begin{align*}
W_{\mathrm{loc}}^{ss}(x) &= \{ y\in U\,|\, \lim_{n\to\infty}\|F^n(y)-F^n(x)\|=0 \},\\
W_{\mathrm{loc}}^{uu}(x) &= \{ y\in U\,|\, \lim_{n\to-\infty}\|F^n(y)-F^n(x)\|=0\},
\end{align*}
where $U$ is a neighborhood of $M$ in $\mathbb{R}^N$. These leaves are $C^r$ submanifolds and tangent to $N_x^s$, $N_x^u$, respectively, at $x\in M$. 
\end{theorem}

\section{Modification of boundaries}
\label{sec_app:boundary}
\setcounter{equation}{0}
Theorems such as provided in~\cite{Hirsch:77:IM} are by default not suitable for manifolds with non-empty boundary. 
In the following, the concept of overflowing and inflowing invariant manifolds of Fenichel is employed to show the existence of the stable and unstable manifolds.
As the name suggests, these overflowing and inflowing manifolds are required to have a strictly outward and respectively inward oriented vector field at the boundary.
As the boundary of $\mathscr{M}(\delta)$ is given by the two periodic orbits of Hamiltonian value $\pm \delta$, $\mathscr{M}(\delta)$ is neither outflowing nor inflowing.
Therefore, the proof is split into four steps.
Firstly, the Hamiltonian vector field $X_H$ is modified near the boundary $\partial \mathscr{M}(\delta)$ to be transverse as proposed in Fenichel~\cite{Wiggins1988}. 
Thereupon, it is shown in step two that for the modified vector field, there exists a hyperbolic splitting of the tangent bundle $T_{\mathscr{M}(\delta)} \M$.
In the third step, it is concluded that there are stable and unstable manifolds of $\mathscr{M}(\delta)$ for the modified vector field $\tilde{X}_H$.
Finally, it is concluded that these manifolds persists under perturbations, i.e. they exist also for $X_H$, which is $C^0$-close to $\tilde{X}_H$.
\subsection{Proof using Boundary Modification}
Let $0<\delta_1<\delta_2<\delta_3<\delta_4$ be sufficiently small such that $\mathscr{M}(\delta_i)$ is well-defined for $i\in \{1,\dots, 4\}$.
It is clear that $\mathscr{M}(0)\subset \mathscr{M}(\delta_1)\subset\dots\subset\mathscr{M}(\delta_4)$ and
\[\mathscr{M}(\delta_i)=\text{int}\mathscr{M}(\delta_i)\cup \partial \mathscr{M}(\delta_i)\subset \mathscr{M}(\delta_{i+1}),\]
where $\text{int}C$ denotes the interior of the set $C$.
There exists $\Psi\in\mathcal{C}^{\infty}(\tangentBundle,\mathbb{R})$ being the sum of two smooth bump functions of $\partial \mathscr{M}(\delta_2)$ and $\partial \mathscr{M}(\delta_3)$ supported in $\text{int}\mathscr{M}(\delta_3)\setminus \mathscr{M}(\delta_1)$ and $\text{int}\mathscr{M}(\delta_4)\setminus \mathscr{M}(\delta_2)$, respectively.
In particular, the function is defined by
\begin{align*}
	\Psi =	\left\{
	\begin{array}{lll}
		0 & \text{ on } \mathscr{M}(\delta_1) \\
		+1 & \text{ on } \partial \mathscr{M}(\delta_2)\\
		-1 & \text{ on } \partial \mathscr{M}(\delta_3) \\
		0 & \text{ on } \tangentBundle\setminus \text{int} \mathscr{M}(\delta_4)
	\end{array},
	\right. 
\end{align*}
see Fig.\ref{fig:bump} for an overview. 
\begin{figure}
	\centering
\begin{tikzpicture}
	\begin{axis}[xmin=-1.5, xmax=5, ymin=-1.2, ymax=1.2,
        xticklabel style={font=\footnotesize},
		xtick = {-1, 0,2,4}, ytick = { -1,0,1}, 
		xticklabels={$\partial\mathscr{M}(\delta_1\!)$, $\partial\mathscr{M}(\delta_2)$, $\partial\mathscr{M}(\delta_3)$,$\partial\mathscr{M}(\delta_4)$}, %
		yscale=0.35, restrict y to domain=-1.2:1.2,
        xscale=1.1,
		axis x line=center, axis y line= left,
		samples=40]
		\addplot[blue, samples=100, smooth, domain=-1.5:-1, very thick]
		plot (\x, { 0 });
		\addplot[blue, samples=100, smooth, domain=-1:1, very thick, label={x}]
		plot (\x, {exp(1-1/(1-x^2)});
		\addplot[blue, samples=100, smooth, domain=-1:1, very thick, label={x}]
		plot (\x+2, {-exp(1-1/(1-x^2)});
		\addplot[blue, very thick, samples=100, smooth, domain=3:4.5]
		plot (\x, {0} );
	\end{axis}
\end{tikzpicture}
	\caption{Function $\Psi$ for the boundary modification.}
	\label{fig:bump}
\end{figure}
Furthermore, the vector field $X_\perp$ denotes
\begin{align*}
	X_\perp = \sum_a \left(
	\partialDer{H}{p_a} \partialDer{}{p_a}
	-\partialDer{H}{x_a} \partialDer{}{x_a}\right) + \partialDer{}{p_0}.
\end{align*}
Near $\mathscr{M}(0)$, it is transverse to the hyperplane defined by a constant Hamiltonian value, which can be easily verified by evaluating $H$ along $X_\perp$, i.e.,
\begin{align*}
	dH(X_\perp) = 1 + f_0(x)-\partialDer{h(x,p)}{p_0}+
	 \left(\partialDer{H}{x_a}\right)^2+\left(\partialDer{H}{p_a}\right)^2,
\end{align*}
where $h(x,p)=\frac{1}{4}G(x,p)^\top R^{-1} G(x,p)$.
Near $\mathscr{M}(0)$, $f_0(x)$ and $\partialDer{h(x,p)}{p_0}$ are zero due to property of the selected transverse coordinates. Furthermore, $\partialDer{H}{p_a}$ and $\partialDer{H}{x_a}$ are zero near $\mathscr{M}(0)$ as $\dot{x}_a =0$ and $\dot{p}_a =0$. Therefore, $dH(X_\perp) >0$ near $\mathscr{M}(0)$ and it follows that $X_\perp$ is an outward pointing vector field for $\mathscr{M}(\delta_2)$ and $\mathscr{M}(\delta_3)$ if $\delta_2$ and $\delta_3$ are sufficiently small.
The modified vector field is defined as
\begin{align*}
	\tilde X_H =  X_H + \Psi X_\perp,
\end{align*}
which alters the Hamiltonian vector field only on $\mathscr{M}(\delta_4)\setminus \mathscr{M}(\delta_1)$ in a smooth way. 
At $\partial \mathscr{M}(\delta_2)$ the Hamiltonian vector field is tangential and $X_\perp$ adds a component to $\tilde X_H$ such that it is an outward pointing vector field for $\mathscr{M}(\delta_2)$. 
Similarly, $\tilde X_H$ is inward pointing at the boundary of $\mathscr{M}(\delta_3)$. 
In the following we are using the notations and conventions of the book \cite{Wiggins:94:NHIMDS}.
For $\delta=\delta_2$ and $\delta=\delta_3$, $\mathscr{M}(\delta)$ is a compact manifold with boundary.
One can now show that the splitting of $\mathscr{M}(\delta)$
\begin{align*}
	T_{\mathscr{M}(\delta)}\M = T\mathscr{M}(\delta) \oplus N^s \oplus N^u
\end{align*}
is hyperbolic by verifying that for $y\in\mathscr{M}(\delta)$
\begin{align*}
	\lambda^u(y) = \inf \left\{ a | \frac{\|u_{-t}\|}{\|u_0\|}/a^t\rightarrow 0, t\rightarrow \infty, \forall u_0 \in N^u_y \right\}=\overline{\lambda}
\end{align*}
and
\begin{align*}
	\nu^s(y) = \inf \left\{ a| \frac{\|w_{0}\|}{\|w_{-t}\|}/a^t\rightarrow 0, t\rightarrow \infty, \forall w_0 \in N^s_y \right\}=\overline{\lambda},
\end{align*} 
for some $\overline{\lambda}<1$ (by shrinking $\epsilon$ if needed).
Similarly, it follows that $\sigma(p)=0$ for all $y\in\mathscr{M}(\delta)$.
$\mathscr{M}(\delta_2)$ is a compact connected $\mathcal{C}^r$ manifold with boundary with hyperbolic splitting and overflowing invariant under $\tilde{X}_H$.
As $\lambda^u<1$, $\nu^s<1$ and $\sigma=0<\frac{1}{r}$ there exists an overflowing invariant manifold $W^u(\mathscr{M}(\delta_2))$ containing $\mathscr{M}(\delta_2)$ by Theorem 1.3.6. of \cite{Wiggins:94:NHIMDS}.
Similarly, $\mathscr{M}(\delta_3)$ is overflowing invariant under the time reversed flow $-\tilde{X}_H$.
Therefore, there is a stable manifold $W^s(\mathscr{M}(\delta_3))$ containing $\mathscr{M}(\delta_3)$.
As $X_H = \tilde{X}_H$ on $\mathscr{M}(\delta_1)$, it follows from Theorem 1.3.6 that there are $W^s(\mathscr{M}(\delta_1))$ and $W^u(\mathscr{M}(\delta_1))$ for the original system with flow $X_H$.

\begin{IEEEbiography}[{\includegraphics[width=1in,height=1.25in,clip,keepaspectratio]{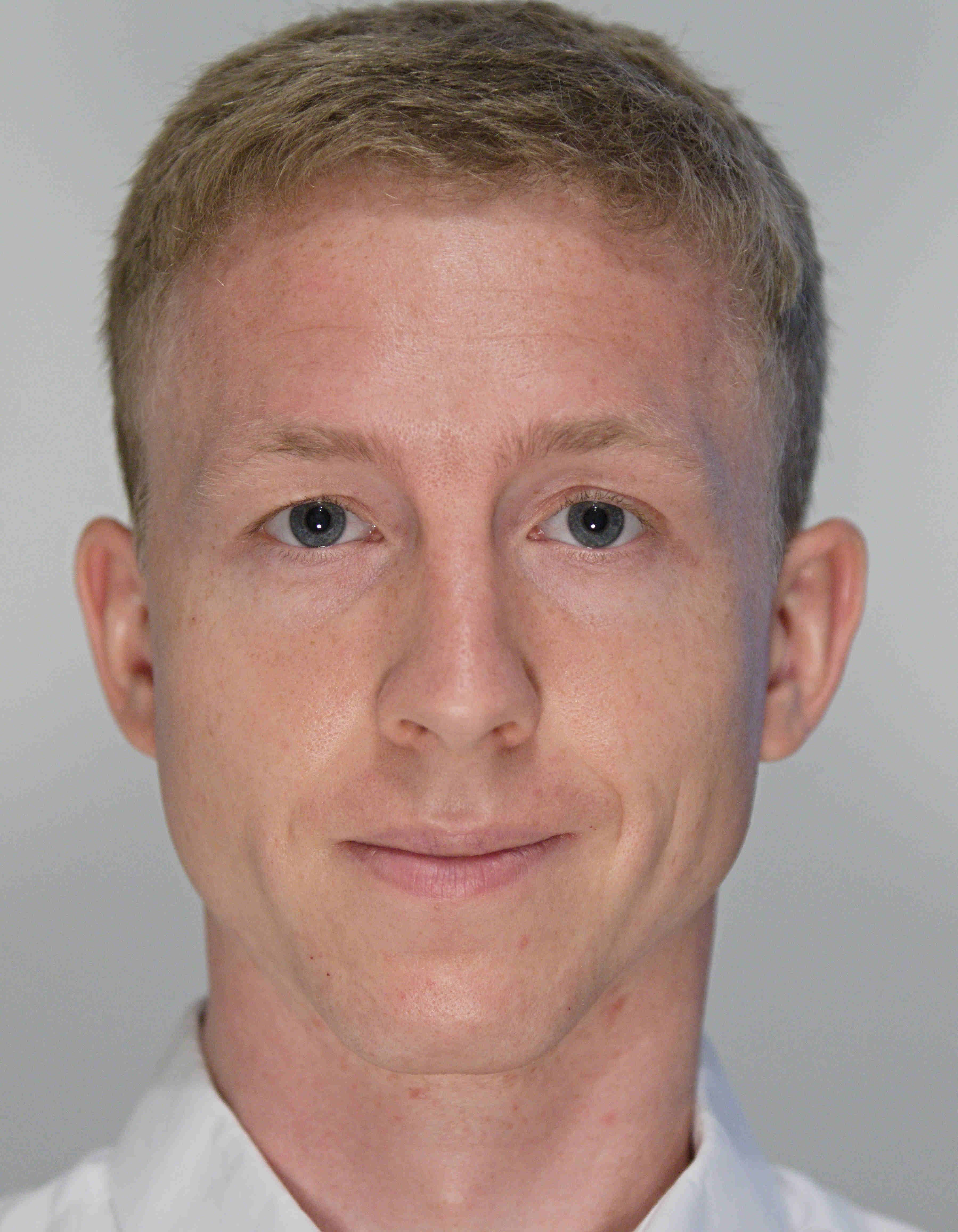}}]{Fabian Beck} received a B.Sc. and M.Sc. degree in Technical Cybernetics and System Theory from the TU Ilmenau in 2016 and 2018, respectively. Since 2018, he has been a research scientist with the Institute of Robotics and Mechatronics of the German Aerospace Center (DLR). Since 2023, he is also with the School of Computation, Information and Technology of the Technical University of Munich (TUM). His research interests include impedance and torque control as well as the control of underactuated systems using methods from geometric analysis and optimal control.
\end{IEEEbiography}

\begin{IEEEbiography}[{\includegraphics[width=1in,height=1.25in,clip,keepaspectratio]{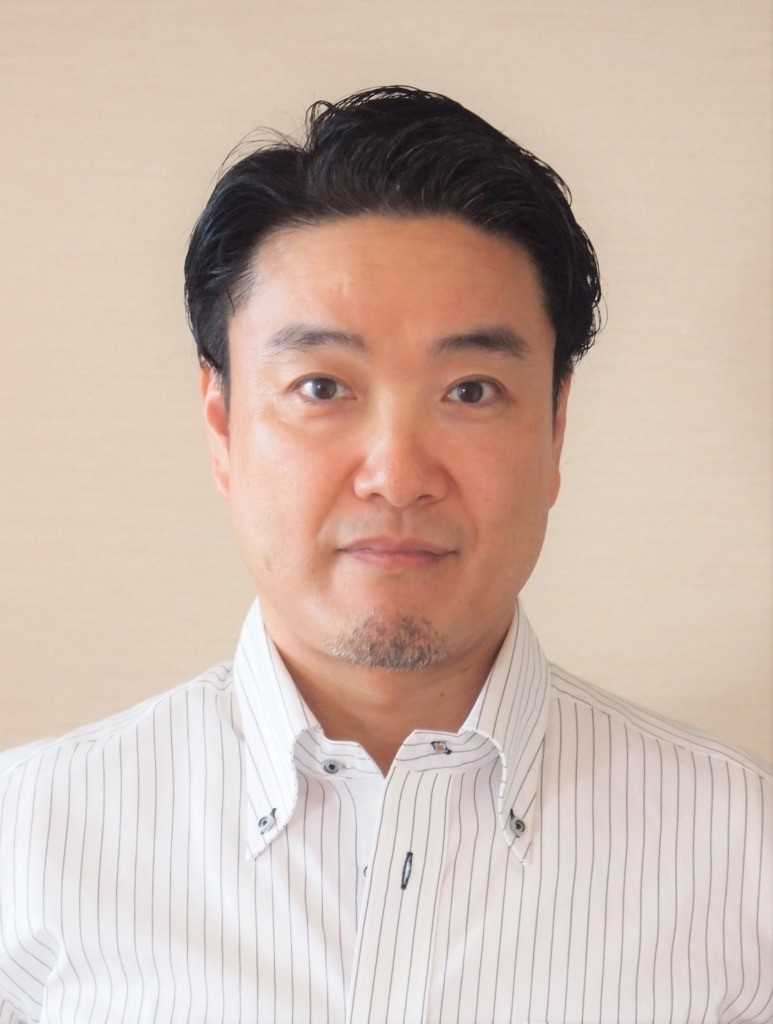}}]{Noboru Sakamoto} received the B.Sc. degree in mathematics from Hokkaido University and M.Sc. and Ph.D. degrees in aerospace engineering from Nagoya University, in 1991, 1993, and 1996, respectively. From 1996 until 2015, he held positions in the Graduate School of Engineering in Nagoya University. Currently, he is a Professor with the Science and Engineering of Nanzan University in Nagoya, Japan. He has held
visiting research positions at University of Groningen, The Netherlands, in 2005 and 2006 and at University of Deusto, Spain, in 2018 and 2019. He received the SICE Best Paper Prizes in 1997, 2006, 2008 and 2011 and Kimura Prize in 2016. His research interests include nonlinear control theory, control of chaotic systems, dynamical system theory and control applications for mechatronics and aerospace engineering.
\end{IEEEbiography}
\end{document}